\documentclass[12pt]{amsart}

\usepackage{latexsym,amscd,amssymb}
\usepackage[dvips]{color,graphicx}
\usepackage{enumerate}
\usepackage[margin=1.5in]{geometry}

\theoremstyle{plain}
  \newtheorem{theorem}{Theorem}[section]
  \newtheorem{proposition}[theorem]{Proposition}
  \newtheorem{lemma}[theorem]{Lemma}
  \newtheorem{corollary}[theorem]{Corollary}
  
\theoremstyle{definition}
  \newtheorem{definition}[theorem]{Definition}
  \newtheorem{example}[theorem]{Example}

 \theoremstyle{remark}
  \newtheorem{remark}[theorem]{Remark}

\numberwithin{equation}{section}

\newcommand{\reals}{{\mathbb R}}
\newcommand{\cplxes}{{\mathbb C}}
\newcommand{\FF}{{\mathbb F}}
\newcommand{\DD}{\Delta}
\newcommand{\WHD}{\widehat{\DD}}
\newcommand{\LK}{\mathrm{link}}
\begin{document}

\title[Shellable complexes and diagonal arrangements]
{Shellable complexes and topology of diagonal arrangements}

\author{Sangwook Kim}

\address{School of Mathematics\\
University of Minnesota\\
Minneapolis, MN 55455, USA}
\curraddr{Department of Mathematical Sciences\\
George Mason University\\
Fairfax, VA 22030, USA}
\email{skim22@gmu.edu}


\thanks{This work forms part of the author's doctoral dissertation
at the University of Minnesota, supervised by Vic Reiner,
and partially supported by NSF grant DMS-0245379.}

\begin{abstract}
   We prove that if a simplicial complex $\DD$ is shellable, 
   then the intersection lattice $L_\DD$ for the corresponding diagonal arrangement
   $\mathcal{A}_\DD$ is homotopy equivalent to a wedge of spheres.
   Furthermore, we describe precisely the spheres in the wedge,
   based on the data of shelling.
   Also, we give some examples of diagonal arrangements $\mathcal{A}$
   where the complement $\mathcal{M}_\mathcal{A}$ is $K(\pi,1)$, 
   coming from rank $3$ matroids.
\end{abstract}

\maketitle

\tableofcontents

\section{Introduction}
\label{sec-introduction}

Consider $\reals^n$ with coordinates $u_1, \ldots, u_n$.
A \emph{diagonal subspace} $U_{i_1 \cdots i_r}$ is a linear subspace of the form 
$u_{i_1} = \cdots = u_{i_r}$.
A \emph{diagonal arrangement} (or a \emph{hypergraph arrangement}) $\mathcal{A}$ is
a finite set of diagonal subspaces of $\reals^n$.
Throughout this paper, we assume that for any $H_1, H_2 \in \mathcal{A}$, $H_1$ is not included in $H_2$.

For a simplicial complex $\DD$ on a vertex set $[n] := \{1, 2, \dots, n \}$ such that $\dim \DD \le n-3$,
one can associate a diagonal arrangement $\mathcal{A}_\DD$ in $\reals^n$ as follows.
For a facet $F$ of $\DD$, let $U_{\overline{F}}$ be the diagonal subspace
$u_{i_1} = \cdots = u_{i_r}$ where $\overline{F} = [n]-F = \{ i_1, \dots, i_r \}$.
Define
$$
\mathcal{A}_\DD = \{ U_{\overline{F}} \mid F \text{ is a facet of } \DD \}.
$$
If every subspace in a diagonal arrangement $\mathcal{A}$ in $\reals^n$ has the form
$U_{i_1 \cdots i_r}$ with $r \ge 2$, one can find a simplicial complex $\DD$ on $[n]$
satisfying $\mathcal{A} = \mathcal{A}_\DD$.

Two important spaces associated with an arrangement $\mathcal{A}$ of linear subspaces
in $\reals^n$ are 
$$
\mathcal{M}_\mathcal{A} = \reals^n - \bigcup_{H \in \mathcal{A}} H
\quad \text{ and } \quad
\mathcal{V}_\mathcal{A}^\circ = \mathbb{S}^{n-1} \cap \bigcup_{H \in \mathcal{A}} H,
$$
called the \emph{complement} and the \emph{singularity link}.

We are interested in the topology of $\mathcal{M}_\mathcal{A}$ and
$\mathcal{V}_\mathcal{A}^\circ$ for a diagonal arrangement
$\mathcal{A}$.
Diagonal arrangements arise in connection with important questions in many different fields.
In computer science,
Bj\"{o}rner, Lov\'{a}sz and Yao~\cite {BjornerLovaszYao} found lower bounds on the complexity
of $k$-equal problems using the topology of diagonal arrangements
(see also \cite{BjornerLovasz}).
In group cohomology,
it is well known that $\mathcal{M}_{\mathcal{B}_n}$ for the braid arrangement 
$\mathcal{B}_n$ in $\cplxes^n$ is a $K(\pi, 1)$ space with the fundamental group 
isomorphic to the pure braid group \cite {FadellNeuwirth}.
Khovanov~\cite {Khovanov} showed that $\mathcal{M}_{\mathcal{A}_{n,3}}$ for 
the $3$-equal arrangement $\mathcal{A}_{n,3}$ in $\reals^n$ is also a $K(\pi, 1)$ space.

Note that $\mathcal{M}_\mathcal{A}$ and $\mathcal{V}_\mathcal{A}^\circ$
are related by Alexander duality as follows:
\begin{eqnarray}
\label{complement-and-link}
   H^i (\mathcal{M}_{\mathcal{A}} ; \FF) 
   = H_{n-2-i} ( \mathcal{V}_\mathcal{A}^\circ  ; \FF)
   \qquad (\FF \text{ is any field})
\end{eqnarray}

In the mid-1980's Goresky and MacPherson~\cite {GoreskyMacpherson} found 
a formula for the Betti numbers of $\mathcal{M}_\mathcal{A}$, while
the homotopy type of $\mathcal{V}_\mathcal{A}^\circ$ was computed by 
Ziegler and \v{Z}ivaljevi\'{c}~\cite {ZieglerZivaljevic}
(see Section~\ref {set-homotopy}).
The answers are phrased in terms of the lower intervals in the \emph{intersection lattice}
$L_{\mathcal{A}}$ of the subspace arrangement $\mathcal{A}$,
that is the collection of all nonempty
intersections of subspaces of $\mathcal{A}$ ordered by reverse inclusion.
For general subspace arrangements, these lower intervals in $L_\mathcal{A}$ can 
have arbitrary homotopy type (see~\cite [Corollary 3.1] {ZieglerZivaljevic}).

Our goal is to find a general sufficient condition for the intersection lattice
$L_{\mathcal{A}}$ of a diagonal arrangement $\mathcal{A}$ to be well behaved.
Bj\"{o}rner and Welker~\cite{BjornerWelker} showed that $L_{\mathcal{A}_{n,k}}$ 
has the homotopy type of a wedge of spheres, where 
$\mathcal{A}_{n,k}$ is the $k$-equal arrangement consisting of all $U_{i_1 \cdots i_k}$
for all $1 \le i_1 < \cdots < i_k \le n$ (see Section~\ref {sec-earlier-results}),
and Bj\"{o}rner and Wachs~\cite{BjornerWachs1} showed that $L_{\mathcal{A}_{n,k}}$
is shellable.
More generally, Kozlov~\cite {Kozlov} showed that 
$L_{\mathcal{A}}$ is shellable if $\mathcal{A}$ 
satisfies certain technical conditions (see Section~\ref {sec-earlier-results}).
Suggested by a homological calculation (Theorem~\ref{nonpure-by-walls} below),
we prove the following main result, capturing the homotopy type assertion from 
\cite {Kozlov} (see Section~\ref{sec-homotopy}).

\begin{theorem}
\label{main-theorem}
   Let $\DD$ be a shellable simplicial complex.
   Then the intersection lattice $L_\DD$ for the diagonal arrangement 
   $\mathcal{A}_\DD$ is homotopy equivalent to a wedge of spheres.
\end{theorem}

Furthermore, one can describe precisely the spheres in the wedge, based on the shelling data.
Let $\DD$ be a simplicial complex on $[n]$ with a shelling order $F_1, \dots, F_q$ on its facets.
Let $\sigma$ be the intersection of all facets, and $\bar{\sigma}$ its complement.
Let $G_1 = F_1$ and for each $i \ge 2$, 
let $G_i$ be the face of $F_i$ obtained by intersecting 
the walls of $F_i$ that lie in the subcomplex generated by $F_1, \dots, F_{i-1}$,
where a \emph{wall of} $F_i$ is a codimension $1$ face of $F_i$.
An (\emph{unordered}) \emph{shelling-trapped decomposition} 
(\emph{of  $\bar{\sigma}$ over $\DD$}) is defined 
to be a family $\{ (\bar{\sigma}_1, F_{i_1}), \dots, (\bar{\sigma}_p, F_{i_p})\}$
such that $\{ \bar{\sigma}_1, \dots, \bar{\sigma}_p \}$ is a 
decomposition of $\bar{\sigma}$ as a disjoint union
$$
\bar{\sigma} = \bigsqcup_{j=1}^p \bar{\sigma}_j
$$
and $F_{i_1} , \dots , F_{i_p}$ are facets of $\DD$
such that $G_{i_j} \subseteq \sigma_j \subseteq F_{i_j}$ for all $j = 1, \dots, p$.
Then the wedge of spheres in Theorem~\ref{main-theorem} 
consists of $(p-1)!$ copies of spheres of dimension
$$
p(2-n) + \sum_{j=1}^p |F_{i_j}| + |\bar{\sigma}|-3
$$
for each shelling-trapped decomposition 
$D = \{ (\bar{\sigma}_1, F_{i_1}), \dots, (\bar{\sigma}_p, F_{i_p})\}$ of $\bar{\sigma}$.
Moreover, for each shelling-trapped decomposition $D$ of $\bar{\sigma}$ and 
a permutation $w$ of $[p-1]$, there exists a saturated chain $C_{D, w}$
(see Section~\ref{sec-chain})
such that removing the simplices corresponding to these chains leaves 
a contractible simplicial complex.

The following example shows that the intersection lattice in Theorem \ref{main-theorem} 
is not shellable in general, even though it has the homotopy type of a wedge of spheres.

\begin{example}
\label{ex-not-shellable}
   Let $\DD$ be a simplicial complex on $[8]$ with a shelling
   $123456$ (short for $\{ 1, 2, 3, 4, 5, 6 \}$), $127, 237, 137, 458, 568, 468$.
   Then the order complex of the upper interval $(U_{78}, \hat{1})$ is a disjoint union
   of two circles, hence is not shellable. 
   Therefore, the intersection lattice $L_\DD$ for the diagonal arrangement
   $\mathcal{A}_\DD$ is also not shellable.
   The intersection lattice $L_\DD$ is shown in Figure \ref{nonshellable}
   (thick lines represent the open interval $(U_{78}, \hat{1})$).
   In Figures, the subspace $U_{i_1 \dots i_r}$ is labeled by $i_1 \dots i_r$.
   Also note that a facet $F$ of $\DD$ corresponds to the subspace $U_{[n] - F}$. 
   For example, the facet $127$ corresponds to $U_{34568}$.

   \begin{figure}
   \begin{center}
   \includegraphics[width=0.9\textwidth]{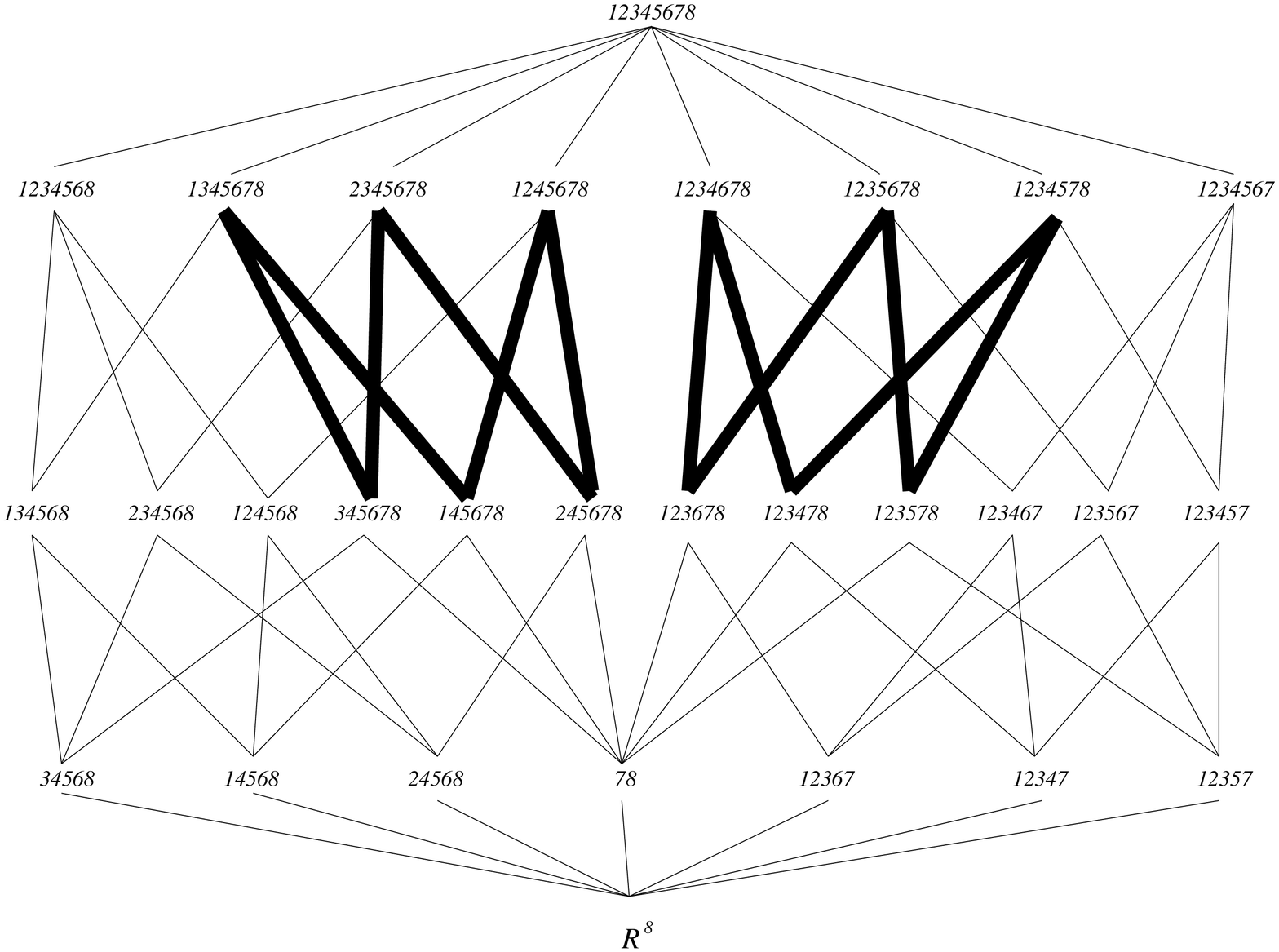}
   \begin{caption}
      {The intersection lattice $L_\DD$ for Example~\ref{ex-not-shellable}
      ($U_{i_1 \dots i_r}$ is labeled by $i_1 \dots i_r$.)}
      \label{nonshellable}
   \end{caption}
   \end{center}
   \end{figure}
\end{example}

The next example shows that there are nonshellable simplicial complexes 
whose intersection lattices are shellable.

\begin{example}
   Let $\DD$ be a simplicial complex on $[4]$ whose facets are $12$ and $34$.
   Then $\DD$ is not shellable. But the order complex of $\overline{L}_\DD$ consists
   of two vertices, hence is shellable.
\end{example}

In Section~\ref{sec-earlier-results}, we give Kozlov's result and show that 
its homotopy type consequence is a special case of Theorem~\ref{main-theorem}.
Also, we give a new proof of Bj\"{o}rner and Welker's result 
using Theorem~\ref{main-theorem}.
In Section~\ref {sec-homotopy}, we prove Theorem~\ref{main-theorem}.
In Section~\ref {set-homotopy}, we deduce from it the homotopy type and 
the homology of the singularity link 
(and hence the homology of the complement)
of a diagonal arrangement
$\mathcal{A}_{\DD}$ for a shellable simplicial complex $\DD$.
In Section~\ref{sec-examples}, we give some examples in which 
$\mathcal{M}_\mathcal{A}$ are $K(\pi,1)$, coming from matroids.

\section{Basic notions and definitions}
\label{sec-basic-notions}

In this section, we give a few definitions used throughout this paper.

First, we start with the definition of (nonpure) shellability of simplicial complexes;
see \cite{BjornerWachs1}, \cite{BjornerWachs2} for further background on this notion.
\begin{definition}
   A simplicial complex is \emph{shellable} if its facets can be arranged in linear order
   $F_1, F_2, \ldots, F_q$ in such a way that % the subcomplex
   $(\bigcup_{i = 1}^{k-1} 2^{F_i}) \cap 2^{F_k}$
   is pure and $(\dim F_k -1)$-dimensional for all $k = 2, \ldots, q$, where
   $2^F = \{ G |  G \subseteq F \}$.
   Such an ordering of facets is called a \emph{shelling order} or \emph{shelling}.
\end{definition}

There are several equivalent definitions of shellability.
The following restatement of shellability is often useful.

\begin{lemma}{\cite [Lemma 2.3] {BjornerWachs1}}
\label{restatement-of-shellability}
   A linear order $F_1, F_2, \ldots, F_q$ of facets of a simplicial complex is a shelling
   if and only if for every $F_i$ and $F_k$ with $F_i < F_k$, 
   there is a facet $F_j <F_k$ such that
   $F_i \cap F_k \subseteq F_j \cap F_k \lessdot F_k$,
   where $G \lessdot F$ means that $G$ has codimension $1$
   in $F$.
\end{lemma}

It is well known that a pure $d$-dimensional shellable simplicial complex has
the homotopy type of a wedge of $d$-spheres.
Bj\"{o}rner and Wachs~\cite {BjornerWachs1} generalized this result to the nonpure case,
i.e., a nonpure shellable simplicial complex has the homotopy type of a wedge of spheres.
However, these spheres need not be equidimensional.

The \emph{link} of a face $\sigma$ of a simplicial complex $\DD$ is
$$
\LK_{\DD} \sigma = \{ \tau \in \DD \mid \tau \cap \sigma = \emptyset
\text{ and } \tau \cup \sigma \in \DD \}.
$$
Bj\"{o}rner and Wachs~\cite{BjornerWachs2} showed that shellability is inherited
by all links of faces in a simplicial complex.

\begin{proposition}
\label{link-shellable}
   If $\DD$ is shellable, then so is $\LK_\DD \sigma$ for all faces $\sigma \in \DD$,
   using the induced order on facets of $\LK_\DD \sigma$.
\end{proposition}

Now we give the definition of the order complex of a poset which we will use
frequently.

\begin{definition}
   The \emph{order complex} $\DD (P)$ of a poset $P$ is the simplicial complex
   whose vertices are the elements of $P$ and
   whose faces are the chains of $P$.
\end{definition}

For the order complex $\DD ((x, y))$ of an open interval $(x, y)$,
we will use the notation $\DD(x,y)$.
When we say that a finite lattice $L$ with bottom element $\hat{0}$
and top element $\hat{1}$ has some topological properties, such as purity,
shellability and homotopy type, it means the order complex of 
$\overline{L} := L - \{ \hat{0}, \hat{1} \}$ 
has those properties.

\section{Special cases that were known}
\label {sec-earlier-results}

In this section, we give Kozlov's theorem and show how its consequence for
homotopy type follows from Theorem~\ref{main-theorem}.
Also, we give a new proof of Bj\"{o}rner and Welker's theorem about
the intersection lattice of the $k$-equal arrangements
using Theorem~\ref{main-theorem}.

Kozlov~\cite {Kozlov} showed that $\mathcal{A}_{\DD}$ has
shellable intersection lattice if $\DD$ satisfies some conditions.
This class includes $k$-equal arrangements and all other diagonal arrangements
for which the intersection lattice was proved shellable.

\begin{theorem}{\cite [Corollary 3.2]{Kozlov}}
\label {Kozlov-corollary}
   Consider a partition of 
   $
   [n] = E_1 \sqcup \cdots \sqcup E_r
   $
   such that
   $\max E_i < \min E_{i+1}$ for $i = 1, \ldots, r-1$.
   Let 
   $$
   f : \{ 1, 2, \ldots, r \} \to \{ 2, 3, \ldots \}
   $$
   be a nondecreasing map.
   Let $\DD$ be a simplicial complex on $[n]$ such that
   $F$ is a facet of $\DD$ if and only if
   \begin{enumerate}
      \item
      $| E_i - F | \le 1$ for $i = 1, \ldots, r$;
      \item
      if $\min \overline{F} \in E_i$, then $| F | = n - f(i)$.
   \end{enumerate}
   Then the intersection lattice for $\mathcal{A}_{\DD}$ is shellable.
   In particular, $L_{\DD}$ has the homotopy type of a wedge of spheres.
\end{theorem}

\begin{proposition}
   $\DD$ in Theorem~\ref {Kozlov-corollary} is shellable.
\end{proposition}

\begin{proof}
   We claim that a shelling order is $F_1, F_2, \ldots, F_q$ such that the words
   $w_1, w_2, \ldots, w_q$ are in lexicographic order, where $w_i$ is the
   increasing array of elements in $\overline{F}_i$.
   Let $F_s, F_t$ be two facets of $\DD$ with $1 \le s < t \le q$.
   Then $w_s \prec_{lex} w_t$.
   Let $m$ be the first number in $[r]$ such that $E_m - F_s \ne E_m - F_t$.
   Construct the word $w$ as follows:
   \begin{enumerate}
      \item
      $w \cap E_i = w_s \cap E_i$ for $i = 1, \ldots, m$;
      \item
      for $i = m+1, \ldots , q$,
      $$
      w \cap E_i = \left\{ 
      \begin{array}{cl}
         w_t \cap E_i & \text{ if } |w \cap (\cup_{j=1}^i E_j)| \le f(l);\\
         \emptyset & \text{ otherwise},
      \end{array}
      \right.
      $$
      where $\min w_s \in E_l$.
   \end{enumerate}
   Note that the length of $w$ is $f(l)$ and $w \prec_{lex} w_t$.
   Let $F $ be the set of all elements which do not appear in $w$. 
   Since $F$ satisfies the two conditions from Theorem~\ref{Kozlov-corollary},
   $F$ is a facet of $\DD$.
   Since $F \cap F_t = F_t - (E_m - F_s)$ and $E_m - F_s$ is a subset of $F_t$ of 
   size $1$, $F \cap F_t$ has codimension $1$ in $F_t$.
   Also $F_s \cap F_t \subseteq F \cap F_t$.
   Hence $F_1, F_2, \ldots, F_q$ is a shelling by Lemma~\ref{restatement-of-shellability}.
\end{proof}

\begin{example}
\label{Kozlov-example}
   Consider the partition of 
   $$
   [7] = \{ 1 \} \sqcup \{ 2, 3 \} \sqcup \{ 4 \} \sqcup \{ 5, 6, 7 \}
   $$
   and the function $f$ given by $f(1) = 2, f(2) = 3, f(3) = 4$, and $f(4) = 5$.
   Then the facets of the simplicial complex $\DD$ that satisfy the conditions from
   Corollary~\ref {Kozlov-corollary} and the corresponding words can be found in 
   Table~\ref {table-Kozlov-example}.
   \begin{table}
   \begin{tabular}{|c|c|c||c|c|c||c|c|c|}
      \hline
      $\min \overline{F}$    & $F$  & $w$ & $\min \overline{F}$    & $F$  & 
      $w$ & $\min \overline{F}$    & $F$  & $w$\\
      \hline
      $1$ & $23456$ & $17$ & $2$ & $1356$ & $247$ & $3$ & $1256$ & $347$ \\
              & $23457$ & $16$ &        & $1357$ & $246$ &          & $1257$ & $346$ \\
              & $23467$ & $15$ &        & $1367$ & $245$ &          & $1267$ & $345$ \\
              & $23567$ & $14$ & & & & & &\\
              & $24567$ & $13$ & & & & & &\\
              & $34567$ & $12$ & & & & & &\\
      \hline
   \end{tabular}
   \vspace{2mm}
      \begin{caption}
      {Facets of $\DD$ from Example~\ref{Kozlov-example} and corresponding words}
      \label{table-Kozlov-example}
   \end{caption}
   \end{table}
   Thus the ordering $34567,$ $24567, $ $23567, $ $23467,$ $23457,$ 
   $23456,$ $1367,$ $1357,$ $1356,$ $1267,$ $1257$
   and $1256$ is a shelling for $\DD$.
\end{example}

One can also use Theorem~\ref{main-theorem} to recover the following theorem of 
Bj\"{o}rner and Welker~\cite{BjornerWelker}.

\begin{theorem}
   The order complex of the intersection lattice $L_{\mathcal{A}_{n,k}}$ has
   the homotopy type of a wedge of spheres consisting of
   $$
   (p-1)! \sum_{0 = i_0 \le i_1 \le \cdots \le i_p = n - pk} \,\,\,
   \prod_{j=0}^{p-1} \binom{n - jk - i_j -1}{k-1} (j+1)^{i_{j+1} - i_j}
   $$
   copies of $(n-3-p(k-2))$-dimensional spheres for $1 \le p \le \lfloor \frac{n}{k} \rfloor$.
\end{theorem}

\begin{proof}
   It is clear that $\mathcal{A}_{n, k} = \mathcal{A}_{\DD_{n, n-k}}$, where
   $\DD_{n, n-k}$ is a simplicial complex on $[n]$ whose facets are 
   all $n-k$ subsets of $[n]$.
   Here, $\sigma = \emptyset$, and so $\bar{\sigma} = [n]$.
   By ordering the elements of each facet in increasing order,
   the lexicographic order of facets of $\DD_{n, n-k}$ gives a shelling.
   Also, one can see that the facets of the form
   $F_i = \{ 1, 2, \dots, m, a_{m+1}, \dots, a_{n-k} \}$,
   where $m+1 <  a_{m+1} < \cdots < a_{n-k}$, have $G_i = \{ 1, \dots, m \}$.
   Thus, $G_i \subseteq \sigma_i \subseteq F_i$ implies
   $\overline{F}_i \subseteq \bar{\sigma}_i \subseteq \overline{G}_i = \{ m+1, \dots, n \}.$
   Note that $\min \bar{\sigma} = \min \overline{F} = m+1$ and $\overline{F}$ 
   has $k$ elements.
   Thus, in any shelling-trapped decomposition $[n] = \sqcup_{j=1}^p \bar{\sigma}_j$,
   one has $p \le \lfloor \frac{n}{k} \rfloor$.
   
   Let $1 \le p \le \lfloor \frac{n}{k} \rfloor$ and 
   $0 = i_0 \le i_1 \le \cdots \le i_p = n - pk$.
   We will construct a shelling-trapped family
   $\{ (\bar{\sigma}_1, F_{i_1}), \dots, (\bar{\sigma}_p, F_{i_p})\}$
   as in Theorem \ref{main-theorem}.
   Since $F_{i_1} < \cdots < F_{i_p}$, we have 
   $\min \bar{\sigma}_1 > \cdots > \min \bar{\sigma}_p$.
   In particular, $1 \in \overline{F}_{i_p} \subseteq \bar{\sigma}_p$.
   Thus there are $\binom{n-1}{k-1}$ ways to pick $\overline{F}_{i_p}$ 
   (equivalently, $F_{i_p}$).
   Now suppose that we have chosen $F_{i_p}, \dots, F_{i_{p-j+1}}$.
   We pick $F_{i_{p-j}}$ so that 
   $\min \overline{F}_{i_{p-j}} = \min \bar{\sigma}_{i_{p-j}}$ is the
   $i_j + 1$st element in $[n] - (F_{i_p} \cup \cdots \cup F_{i_{p-j+1}})$.
   Then we have $\binom {n- jk - i_j -1}{k-1}$ ways to choose $F_{i_{p-j}}$.
   For each $j = 1, \dots, p$, there are $i_j - i_{j-1}$ elements in 
   $
   [n] - (F_{i_p} \cup \cdots \cup F_{i_{p-j+1}})
   $
   which are strictly between $\min \overline{F}_{i_{p-j+1}}$ and
   $\min \overline{F}_{i_{p-j}}$ and they must be contained in one of
   $\bar{\sigma}_p, \dots, \bar{\sigma}_{p-j+1}$
   (i.e., there are $j^{i_j - i_{j-1}}$ choices).
   Therefore there are
   $$
   \prod_{j=0}^{p-1} \binom {n- jk - i_j -1}{k-1}
   \prod_{j=1}^p j^{i_j - i_{j-1}}
   = \prod_{j=0}^{p-1} \binom{n - jk - i_j -1}{k-1} (j+1)^{i_{j+1} - i_j}
   $$
   shelling-trapped families.
   By Theorem~\ref{main-theorem}, each of those families contributes $(p-1)!$ copies of 
   spheres of dimension 
   $$
   p(2-n) + \sum_{j=1}^p (n-k) + n - 3 = n-3 - p (k-2).
   $$
\end{proof}

\section{Proof of main theorem}
\label{sec-homotopy}

Theorem~\ref{main-theorem} will be deduced from a more general statement
about homotopy types of lower intervals in $L_{\DD}$,
Theorem~\ref{one-set-nonpure} below.
Throughout this section, we assume that $\DD$ is a simplicial complex 
on $[n]$ with $\dim \DD \le n-3$.

\begin{theorem}
\label{one-set-nonpure}
   Let $F_1, \dots, F_q$ be a shelling of $\DD$
   and $U_{\bar{\sigma}}$ a subspace in $L_{\DD}$ for some subset
   $\bar{\sigma}$ of $[n]$.
   Then $\DD ( \hat{0}, U_{\bar{\sigma}} )$ is homotopy equivalent to 
   a wedge of spheres, consisting of $(p-1)!$ copies of 
   spheres of dimension
   $$
   \delta(D) := p(2-n) + \sum_{j=1}^p |F_{i_j}|+|\bar{\sigma}|-3
   $$
   for each shelling-trapped decomposition 
   $D = \{(\bar{\sigma}_1, F_{i_1}), \dots, (\bar{\sigma}_p , F_{i_p} )\}$
   of $\bar{\sigma}$.
   
   Moreover, for each such shelling-trapped decomposition $D$ and 
   each permutation $w$ of $[p-1]$, one can construct a saturated chain 
   $C_{D, w}$ (see Section~\ref{sec-chain} below), such that
   if one removes the corresponding $\delta(D)$-dimensional simplices for all
   pairs $(D, w)$, the remaining simplicial complex 
   $\WHD (\hat{0}, U_{\bar{\sigma}})$ is contractible.
\end{theorem}

To prove this result, we begin with some preparatory lemmas.

First of all, one can characterize exactly which subspaces lie in $L_{\DD}$ 
when $\DD$ is shellable.
Recall that for $\bar{\sigma} = \{ i_1, \dots, i_r \} \subseteq [n]$, we denote by
$U_{\bar{\sigma}}$ the diagonal subspace of the form $ u_{i_1} = \cdots = u_{i_r}$.
We also use the notation $U_{\bar{\sigma}_1 / \cdots / \bar{\sigma}_p}$ to denote
$U_{\bar{\sigma}_1} \cap \cdots \cap U_{\bar{\sigma}_p}$ for pairwise disjoint subsets 
$\bar{\sigma}_1, \dots, \bar{\sigma}_p$ of $[n]$.

A simplicial complex is called \emph{gallery-connected} if any pair $F, F'$ of facets
are connected by a path
$$
F = F_0, F_1, \dots, F_{r-1}, F_r = F'
$$
of facets such that $F_i \cap F_{i-1}$ has codimension $1$ in $F_i$ for $i = 1, \dots, r$.
Since it is known that Cohen-Macaulay simplicial complexes are gallery-connected,
shellable simplicial complexes are gallery-connected.

\begin{lemma}
\label {element-of-lattice}
   \begin{enumerate}
      \item
      Given any simplicial complex $\DD$ on $[n]$,
      every subspace $H$ in $L_{\DD}$ has the form
      $$
      H = U_{\bar{\sigma}_1 / \cdots / \bar{\sigma}_p}
      $$
      for pairwise disjoint subsets $\bar{\sigma}_1, \ldots, \bar{\sigma}_p$ of $[n]$
      such that $\sigma_i$ can be expressed as an intersection of
      facets of $\DD$ for $i = 1, 2, \dots, p$.
      \item
      Conversely, when $\DD$ is gallery-connected,
      every subspace $H$ of $\reals^n$ that has the above form
      lies in $L_{\DD}$.
   \end{enumerate}
\end{lemma}

\begin{proof}
   To see (1), 
   note that every subspace $H$ in $L_{\DD}$ has the form
   $$
   H = U_{\bar{\sigma}_1 / \cdots / \bar{\sigma}_p}
   $$
   for pairwise disjoint subsets $\bar{\sigma}_1, \ldots, \bar{\sigma}_p$ of $[n]$.
   Since $H = \bigcap_{F \in \mathcal{F}} U_{\overline{F}}$ for some family $\mathcal{F}$
   of facets of $\DD$, 
   $$
   U_{\bar{\sigma}_j} = \bigcap_{F \in \mathcal{F}_j} U_{\overline{F}}
   $$
   for some subfamily $\mathcal{F}_j$ of $\mathcal{F}$ for all $j = 1, \dots, p$.
   Therefore
   $$
   \sigma_j = \bigcap_{F \in \mathcal{F}_j} F
   $$
   for $j = 1, \dots, p$.
      
   For (2), suppose that $H$ has the form
   $$
   H = U_{\bar{\sigma}_1 / \cdots / \bar{\sigma}_p}
   $$
   for pairwise disjoint subsets $\bar{\sigma}_1, \ldots, \bar{\sigma}_p$ of $[n]$
   such that $\sigma_i$ can be expressed as an intersection of
   facets of $\DD$ for $i = 1, 2, \dots, p$.
   It is enough to show the case when $H = U_{\bar{\sigma}}$.
   Since gallery-connectedness is inherited by all links of faces in a simplicial complex, 
   we may assume $\sigma = \bigcap_{F \in \mathcal{F}} F$, 
   where $\mathcal{F}$ is the set of all facets of $\DD$,
   without loss of generality.
   Then $\bar{\sigma} = \bigcup_{F \in \mathcal{F}} \overline{F}$.
   
   We claim that the simplicial complex $\Gamma$ whose facets are
   $\{ \overline{F} \mid F \in \mathcal{F} \}$ is connected.
   Since $\dim \DD \le n-3$, every facet of $\Gamma$ has at least two elements.
   Let $\overline{F}, \overline{F'}$ be two facets of $\Gamma$ with $F < F'$.
   Since $\DD$ is gallery-connected, there is $F=F_1, F_2, \dots, F_k = F'$ such that
   $F_i \cap F_{i-1}$ has codimension $1$ in $F_i$ for all $i = 2, \dots, k$.
   Thus $\overline{F}_i$ and $\overline{F}_{i-1}$ share at least one vertex for all
   $i = 2, \dots, k$.
   This implies that $\overline{F}$ and $\overline{F'}$ are connected.
   Hence $\Gamma$ is connected.
   
   Therefore $U_{\bar{\sigma}} = \bigcap_{F \in \mathcal{F}} U_{\overline{F}}$.
\end{proof}

The next example shows that the conclusion of Lemma~\ref{element-of-lattice}(2) 
can fail when $\DD$ is not assumed to be gallery-connected.

\begin{example}
   Let $\DD$ be a simplicial complex with two facets $123$ and $345$.
   Then $\DD$ is not gallery-connected.
   Since $L_{\DD}$ has only four subspaces $\reals^5, U_{12},$ $U_{45}$ and $U_{12 / 45}$,
   it does not have the subspace $U_{1245}$, even though $\overline{1245} = 3$ is
   an intersection of facets $123$ and $345$ of $\DD$.
   Thus the conclusion of Lemma~\ref{element-of-lattice}(2) fails for $\DD$.
\end{example}

The following easy lemma, whose obvious proof is omitted, 
shows that every lower interval $[\hat{0}, H]$ can be written as
a product of lower intervals of the form $[\hat{0}, U_{\bar{\sigma}}]$.

\begin{lemma}
\label{more-than-one}
   Let $H \in L_{\DD}$ be a subspace of the form
   $$
   H = U_{\bar{\sigma}_1 / \cdots / \bar{\sigma}_p}
   $$
   for pairwise disjoint subsets $\bar{\sigma}_1, \dots, \bar{\sigma}_p$ of $[n]$.
   Then 
   $$
   [\hat{0}, H] = [\hat{0}, U_{\bar{\sigma}_1}] \times \cdots 
   \times [\hat{0}, U_{\bar{\sigma}_p}].
   $$
   In particular,
   $$
   \DD ( \hat{0}, H) = \DD(\hat{0}, U_{\bar{\sigma}_1}) * \cdots *
   \DD(\hat{0}, U_{\bar{\sigma}_p}) * \mathbb{S}^{p-2},
   $$
   where $``*"$ denotes join of spaces.
\end{lemma}

Note that the join of a space $X$ with an empty set equals $X$.

The next lemma, whose proof is completely straightforward and omitted, 
shows that the lower interval $[\hat{0}, U_{\bar{\sigma}}]$ is isomorphic
to the intersection lattice for the diagonal arrangement corresponding to
$\LK_\DD \sigma$.

\begin{lemma}
\label{iso-to-links}
   Let $U_{\bar{\sigma}}$ be a subspace in $L_{\DD}$ for some
   face $\sigma$ of $\DD$.
   Then the lower interval $[\hat{0}, U_{\bar{\sigma}}]$ is isomorphic to the
   intersection lattice of the diagonal arrangement
   $\mathcal{A}_{\LK_{\DD} (\sigma)}$ corresponding to $\LK_{\DD}
   (\sigma)$ on the vertex set $\bar{\sigma}$.
\end{lemma}

The following lemma shows that upper intervals in $L_\DD$ are at least still 
\emph{homotopy equivalent} to the intersection lattice of a diagonal arrangement.

\begin{lemma}
\label{lemma-upper-interval}
   Let $U_{\bar{\sigma}}$ be a subspace in $L_{\DD}$ for some
   face $\sigma$ of $\DD$.
   Then the upper interval $[U_{\bar{\sigma}}, \hat{1}]$ is homotopy equivalent to the
   intersection lattice of the diagonal arrangement
   $\mathcal{A}_{\DD_\sigma}$ corresponding to the simplicial complex $\DD_\sigma$
   on the vertex set $\sigma \cup \{ v \}$ whose facets are obtained 
   in the following ways:
   \begin{enumerate}[\qquad (A)]
      \item
      \label{intersection}
      If $F \cap \sigma$ is maximal among 
      $$
      \{ F \cap \sigma \mid F \text{ is a facet of } \DD \text{ with }
      \sigma \nsubseteq F \text{ and } F \cup \sigma \ne [n] \},
      $$
      then $\widetilde{F} = F \cap \sigma$ is a facet of $\DD_\sigma$.
      \item
      \label{union}
      If a facet $F$ of $\DD$ satisfies $F \cup \sigma = [n]$,
      then $\widetilde{F} = (F \cap \sigma) \cup \{ v \}$ is a facet of $\DD_\sigma$.
   \end{enumerate}
\end{lemma}

\begin{proof}
   We apply a standard crosscut/closure lemma (\cite [Theorem 10.8] {Bjorner95}) 
   saying that 
   a finite lattice $L$ is homotopy equivalent to the sublattice consisting of the joins
   of subsets of its atoms.
   By the closure relation $\psi_1$ on $[U_{\bar{\sigma}}, \hat{1}]$ which sends a subspace to 
   the intersection of all subspaces that lie weakly below it and cover $U_{\bar{\sigma}}$ 
   in $[U_{\bar{\sigma}}, \hat{1}]$,
   one can see that $[U_{\bar{\sigma}}, \hat{1}]$ is homotopy equivalent to 
   the sublattice $L_\sigma$
   generated by the subspaces of $[U_{\bar{\sigma}}, \hat{1}]$ that cover $U_{\bar{\sigma}}$.
   Using the map $\psi_2$ defined by
   $$
   \psi_2(U_{\bar{\tau}}) = \left\{ 
   \begin{array}{cl}
   U_{(\bar{\tau} - \bar{\sigma}) \cup \{ v \}}
   &\text{if } \bar{\sigma} \cap \bar{\tau} \ne \emptyset, \\
   U_{\bar{\tau}}
   &\text{otherwise,}
   \end{array}
   \right.
   $$
   one can see that
   $L_\sigma$ is isomorphic to the intersection lattice $L_{\DD_\sigma}$
   for a simplicial complex $\DD_\sigma$ on the vertex set $\sigma \cup \{ v \}$.
   The facets of $\DD_\sigma$ correspond to the subspaces that cover $U_{\bar{\sigma}}$,
   giving the claimed characterization of facets of $\DD_\sigma$.
   Therefore the map $\psi := \psi_2 \circ \psi_1$ gives the homotopy equivalence between
   $[U_{\bar{\sigma}}, \hat{1}]$ and $L_{\DD_\sigma}$.
\end{proof}

Facets of $\DD_\sigma$ that do not contain $v$ are called \emph{facets of type (\ref{intersection})},
and facets of $\DD_\sigma$ containing $v$ are called \emph{facets of type (\ref{union})}.

\begin{example}
\label{ex-upper-iso}
   Let $\DD$ be a simplicial complex on $[5]$ with facets
   $123,$ $234,$ $35,$ $45$ and let 
   $\sigma = \{ 1, 2, 3 \}$.
   Then $\DD_\sigma$ is a simplicial complex on $\{ 1, 2, 3, v \}$ and 
   its facets are $23$ and $v$. 
   The intersection lattices $L_{\DD}$ and $L_{\DD_\sigma}$ are shown in 
   Figure~\ref{pic-iso-intervals} and it is easy to see that 
   the order complex for $\overline{L}_{\DD_\sigma}$ is homotopy equivalent to 
   the order complex for the open interval $(U_{45}, \hat{1})$ in $L_\DD$.
   Note that the thick lines in Figure~\ref{pic-iso-intervals} (a)
   represent the closed interval $[U_{45}, \hat{1}]$ in $L_\DD$.

   \begin{figure}[!t]
   \begin{center}
      \begin{tabular}{ccc}
         \includegraphics[width=0.5\textwidth]{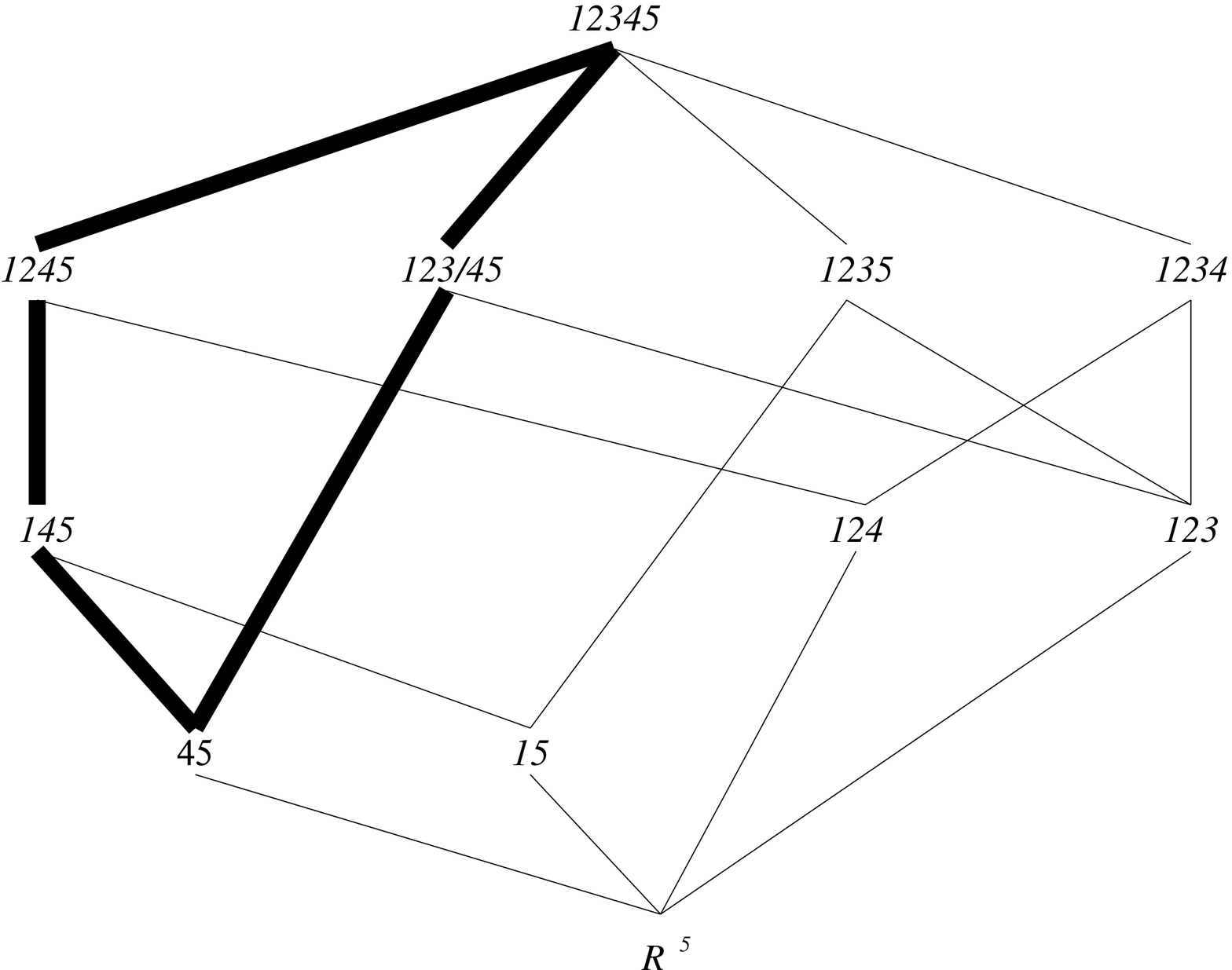}& &
         \includegraphics[width=0.23\textwidth]{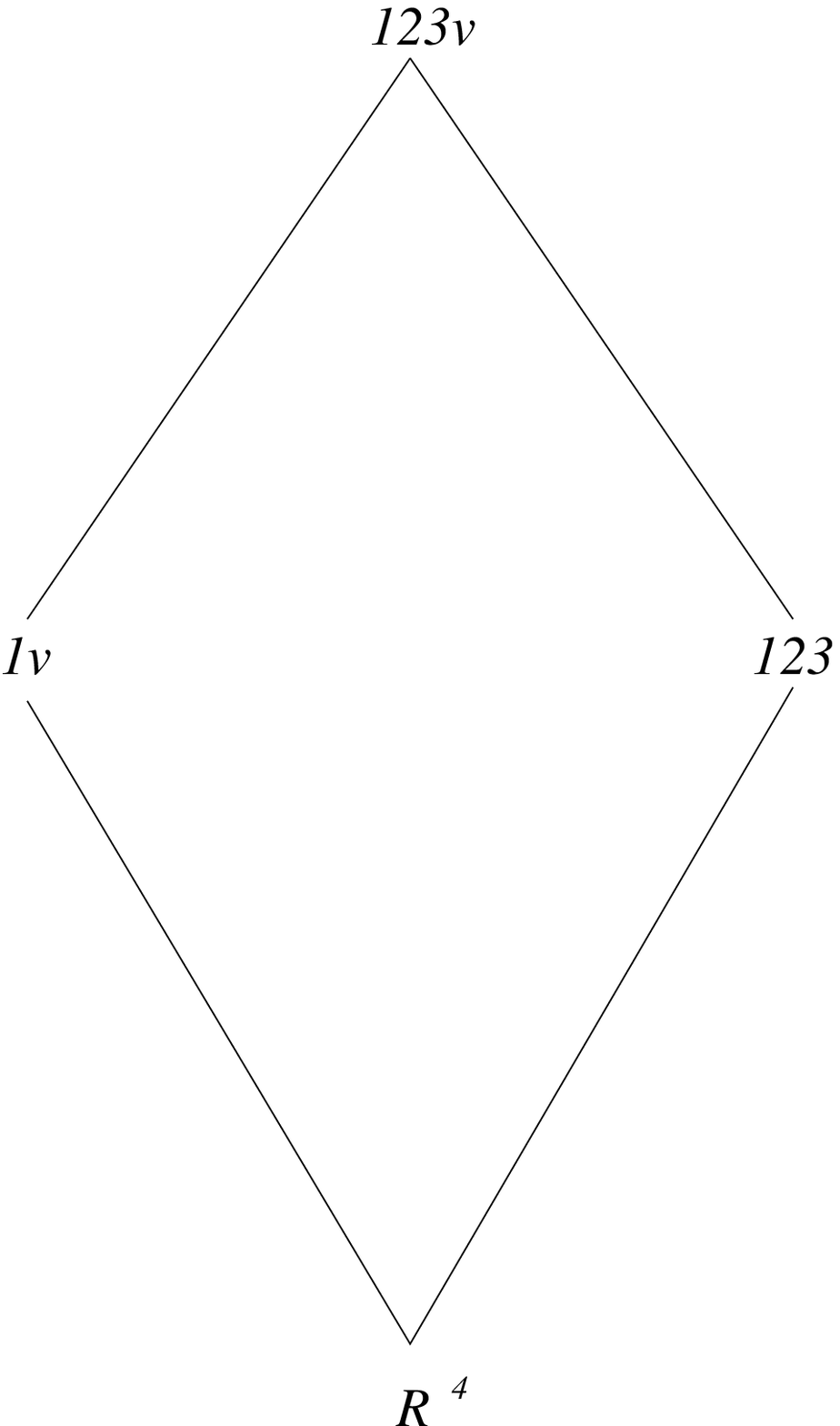}\\
         {\small (a) $L_\DD$} & \hspace{2cm} & {\small (b) $L_{\DD_\sigma}$} \\
      \end{tabular}
      \begin{caption}
         {The intersection lattices for $\DD$ and $\DD_\sigma$}
         \label{pic-iso-intervals}
      \end{caption}
   \end{center}
   \end{figure}
\end{example}

In general, the simplicial complex $\DD_\sigma$ from Lemma~\ref{lemma-upper-interval}
is not shellable, even though $\DD$ is shellable (see Example~\ref{ex-not-shellable}).
However, the next lemma shows that $\DD_\sigma$ is shellable if $\sigma$ is the last facet 
in the shelling order.

\begin{lemma}
\label{upper-last-facet}
   Let $\DD$ be a shellable simplicial complex.
   If $F$ is the last facet in a shelling order of $\DD$,
   then $\DD_F$ is also shellable.
   Moreover, if $\widetilde{F}_i$ is a facet of $\DD_F$ of type (\ref{union}),
   then $\widetilde{G}_i = G_i \cap F$.
\end{lemma}

\begin{proof}
   One can check that the following gives a shelling order on the facets of $\DD_F$.
   First list the facets of type 
   (\ref {intersection}) according to the order of their earliest corresponding facet of $\DD$, 
   followed by the facets of type (\ref{union}) 
   according to the order of their corresponding facet of $\DD$. 
   
   To see the second assertion, 
   let $\widetilde{F}_i$ be a facet of $\DD_F$ of type (\ref{union}), i.e.,
   $$\widetilde{F}_i = (F_i \cap F) \cup \{ v \}$$ for some facet $F_i$ of $\DD$ 
   such that $F \cup F_i = [n]$.
   Then $\widetilde{G}_i = G_i \cap F$ follows from the observations that 
   $F_i \cap F$ is an old wall of $\widetilde{F}_i$, and
   all other old walls of $\widetilde{F}_i$ are 
   $\widetilde{F}_k \cap \widetilde{F}_i$ for some 
   facets $\widetilde{F}_k$ of $\DD_F$ of type (\ref{union}) such that
   $F_k \cap F_i$ is an old wall of $F_i$.
\end{proof}

\begin{example}
\label{ex-upper-iso-last}
   Let $\DD$ be a shellable simplicial complex on $[7]$
   with a shelling $12367,$ $12346,$ $13467,$ $34567,$ $13457,$ $14567, 12345$
   and let $F = 12345$.
   Then $\DD_F$ is a simplicial complex on $\{ 1, 2, 3, 4, 5, v \}$ and 
   its facets are $123v,$ $1234,$ $134v, 345v, 1345, 145v$.
   Since $1234, 1345$ are facets of $\DD_F$ of type (\ref{intersection})
   and $123v,$ $134v,$ $345v,$ $145v$ are facets of $\DD_F$ of type (\ref{union}),
   the ordering $1234,$ $1345,$ $123v,$ $134v,$ $345v,$ $145v$ is a shelling of $\DD_F$.
\end{example}

We next construct the saturated chains appearing in the statement of 
Theorem \ref{one-set-nonpure}.

\subsection{Constructing the chains $C_{D, w}$}
\label{sec-chain}

Let $\DD$ be a simplicial complex on $[n]$ with a shelling $F_1, F_2, \dots, F_q$ and 
let $U_{\bar{\sigma}}$ be a subspace in $L_\DD$.
Let 
$$
D = \{ ( \bar{\sigma}_1, F_{i_1}), \dots, (\bar{\sigma}_p, F_{i_p}) \}
$$
be a shelling-trapped decomposition of $\bar{\sigma}$ with $i_1 < i_2 < \dots < i_p$,
and let $w$ be a permutation on $[p-1]$.

It is well known that the lattice $\Pi_p$ of partitions of the set $[n]$ ordered by refinement
is homotopy equivalent to a wedge of $(p-1)!$ spheres of dimension $p-3$ and 
there is a saturated chain $C_w$ in $\Pi_p$ for each permutation $w$ of $[p-1]$ such that removing 
$\{ \overline{C}_w = C_w - \{ \hat{0}, \hat{1} \} | w \in \mathfrak{S}_{p-1} \}$
from the order complex of $\overline{\Pi}_p$ gives a contractible subcomplex
(see \cite [Example~2.9] {Bjorner80}).

We construct a chain $C_{D, w}$ in $[\hat{0}, U_{\bar{\sigma}}]$ as follows:
\begin{enumerate}
   \item
   Since $G_{i_j} \subseteq \sigma_j \subseteq F_{i_j}$ for all $j = 1, 2, \dots, p$, 
   Lemma~\ref{element-of-lattice} shows that
   the interval $[U_{\overline{F}_{i_1} / \cdots / \overline{F}_{i_p}}, U_{\bar{\sigma}}]$ 
   contains the subspaces $U_{B_1 / \cdots / B_r}$, where
   $B_m = \cup_{k \in K_m} \bar{\tau}_k$ such that
   \begin{itemize}
      \item
      $(K_1 / K_2 / \cdots / K_r)$ is a partition of $[p]$, and
      \item
      $\overline{F}_{i_k} \subseteq \bar{\tau}_k \subseteq \bar{\sigma}_k$ for all $k = 1, 2, \dots, p$.
   \end{itemize}

   Choose a saturated chain $\widehat{C}_w$ in $[U_{\overline{F}_{i_1} / \cdots / \overline{F}_{i_p}}, U_{\bar{\sigma}}]$
   whose covering relations are one of the following types (see \cite{Kozlov}):
   \begin{enumerate}[(a)]
      \item
      $U_{B_1/ B_2 / \cdots / B_r} \lessdot U_{B_1 \cup B_2 / B_3 / \cdots / B_r}$,  
      \item
      $U_{B_1 / \cdots / B_r} \lessdot U_{B_1 \cup \{ a \} / B_2 / \cdots / B_r}$, where
      $a \in \bar{\sigma}_k - \overline{F}_{i_k}$ such that $\overline{F}_{i_k} \subseteq B_1$,
   \end{enumerate}
   where the covering relations of type (a) correspond to the saturated chain $C_w$ in $\Pi_p$, i.e.,
   the covering relation $U_{B_1/ B_2 / \cdots / B_r} \lessdot U_{B_1 \cup B_2 / B_3 / \cdots / B_r}$ appears
   exactly where $(K_1 / K_2 / \cdots / K_r) \lessdot (K_1 \cup K_2 / K_3 / \cdots / K_r)$ 
   appears in $\Pi_p$ (where $\lessdot$ means the covering relation).
   
   \item
   Define a saturated chain $C_{D, w}$ by
   $$
   \hat{0} \lessdot U_{\overline{F}_{i_p}} \lessdot
   U_{\overline{F}_{i_{p-1}} / \overline{F}_{i_p}} \lessdot \dots \lessdot
   U_{\overline{F}_{i_1} / \cdots / \overline{F}_{i_p}}
   $$
   followed by the chain $\widehat{C}_w$.

\end{enumerate}

Note that the length of the chain 
$\overline{C}_{D, w} = C_{D, w} - \{ \hat{0}, U_{\bar{\sigma}} \}$ is
$$
\begin{aligned}
   l(\overline{C}_{D, w}) 
   &= p + (p-1) + \sum_{j=1}^p \left(|\bar{\sigma}_j| - |\overline{F}_{i_j}| \right) - 2\\
   &= p(2-n) + \sum_{j=1}^p |F_{i_j}| + |\bar{\sigma}| -3.
\end{aligned}
$$

\begin{example}
\label{ex-upper-iso-chain}
   Let $\DD$ be the shellable simplicial complex from Example~\ref{ex-upper-iso-last}.
   Then one can see that 
   $$
   D = \{ (45, F_1 = 12367), (123, F_6 = 14567), (67, F_7 = 12345) \}
   $$
   is a shelling-trapped decomposition of $[7]$.
   Let $w$ be the permutation in $\mathfrak{S}_2$ with $w (1) = 2$ and
   $w (2) = 1$.
   Then the maximal chain $C_w$ in $\Pi_3$ corresponding to $w$ is
   $(1 / 2 / 3) \lessdot (1 / 23) \lessdot (123)$.
   One can choose 
   $$
   \widehat{C}_w = U_{45 / 23 / 67} \lessdot U_{45 / 123 / 67} 
                   \lessdot U_{45 / 12367} \lessdot U_{1234567},
   $$
   where the covering $U_{45 / 123 / 67} \lessdot U_{45 / 12367}$ corresponds to $(1 / 2 / 3) \lessdot (1 / 23)$,
   and $U_{45 / 12367} \lessdot U_{1234567}$ corresponds to $(1 / 23) \lessdot (123)$.
   Thus $C_{D, w}$ is the chain
   $$
   \hat{0} \lessdot U_{67} \lessdot U_{23 / 67} \lessdot
   U_{45 / 23 / 67}
   \lessdot U_{45 / 123 / 67} \lessdot
   U_{45 / 12367} \lessdot U_{1234567}.
   $$
   The upper interval $(U_{67}, \hat{1})$ is shown in 
   Figure~\ref{pic-upper-interval} and 
   the chain $\overline{C}_{D, w}$ is represented by thick lines.

   \begin{figure}
   \begin{center}
   \includegraphics[width=0.9\textwidth]{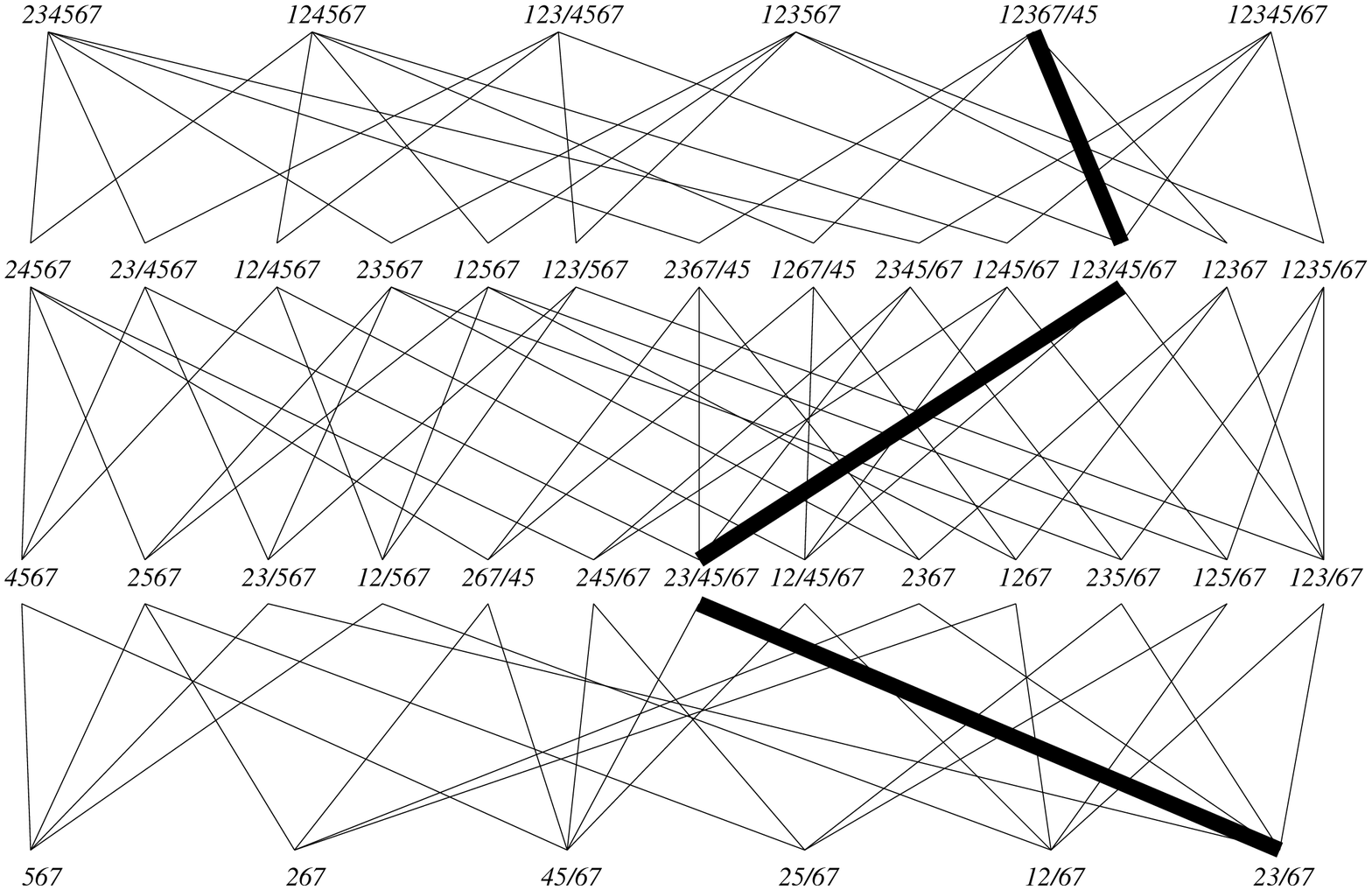}
   \begin{caption}
      {The upper interval $(U_{67}, \hat{1})$ in $L_\DD$}
      \label{pic-upper-interval}
   \end{caption}
   \end{center}
   \end{figure}   
\end{example}

The following lemma gives the relationship between the shelling-trapped decompositions 
of $[n]$ containing $F$ and the shelling-trapped decompositions of $F \cup \{ v \}$.

\begin{lemma}
\label{STD-correspondence}
   Let $\DD$ be a shellable simplicial complex such that the intersection of all facets is empty.
   If $F$ is the last facet in the shelling order of $\DD$,
   then there is a one-to-one correspondence between 
   \begin{itemize}
      \item
      pairs $(D, w)$ of
      shelling-trapped decompositions $D$ of $[n]$ over $\DD$ containing $F$ 
      and $w \in \mathfrak{S}_{|D|-1}$, and 
      \item
      pairs $(\widetilde{D}, \widetilde{w})$ of shelling-trapped decompositions
      $\widetilde{D}$ of $F \cup \{ v \}$ over $\DD_F$ and 
      $\widetilde{w} \in \mathfrak{S}_{|\widetilde{D}|-1}$.
   \end{itemize}
   Moreover, one can choose $\overline{C}_{D, w}$ and $\overline{C}_{\widetilde{D}, \widetilde{w}}$ 
   in such a way that the homotopy equivalence $\psi$ appearing in the proof of Lemma~\ref{lemma-upper-interval} 
   maps the chain $\overline{C}_{D,w} - U_{\overline{F}}$ to the chain $\overline{C}_{\widetilde{D},\widetilde{w}}$. 
\end{lemma}

\begin{proof}
   Let $D = \{ (\bar{\sigma}_1, F_{i_1}), \dots, (\bar{\sigma}_p, F_{i_p}) \}$ be 
   a shelling-trapped decomposition of $[n]$ over $\DD$ with $F_{i_p} = F$
   and let $w$ be a permutation in $\mathfrak{S}_{p-1}$.
   Then $\widetilde{F}_{i_j} = (F_{i_j} \cap F) \cup \{ v \}$ are facets of 
   $\DD_F$ of type (\ref{union}) for $j = 1, \dots, p-1$ and
   $\widetilde{F}_{i_1} < \cdots < \widetilde{F}_{i_{p-1}}$.
   By Lemma~\ref{upper-last-facet}, $\widetilde{G}_{i_j} = G_{i_j} \cap F$ for 
   all $j = 1, \dots, p-1$.
      
   There are two cases to consider:
   
   \vspace{2mm}
   \emph{Case} 1. $\sigma_p \ne F$.
   \vspace{2mm}

   In this case, we will show
   $$
   \widetilde{D} = \{ ( [\bar{\sigma}_p \cap F] \cup \{ v \}, \widetilde{F}),
   ( \bar{\sigma}_1, \widetilde{F}_{i_1}), \dots,
   ( \bar{\sigma}_{p-1}, \widetilde{F}_{i_{p-1}}) \}
   $$
   is a shelling-trapped decomposition of $F \cup \{ v \}$ over $\DD_F$
   ($\widetilde{F}$ will be defined later).

   Define
   $$
   \tilde{\sigma}_j = \left\{ 
   \begin{array}{cl}
      (\sigma_j \cap F) \cup \{ v \} & \text{ for } j = 1, \dots, p-1,\\
      \sigma_j & \text{ for } j = p.
   \end{array}
   \right.
   $$
   For $j = 1, \dots, p-1$,  
   $\widetilde{G}_{i_j} \subseteq \tilde{\sigma}_j \subseteq \widetilde{F}_{i_j}$
   since $G_{i_j} \subseteq \sigma_j \subseteq F_{i_j}$.
   
   Since $G_{i_p} \subseteq \sigma_p$, it must be that $\sigma_p$ is an intersection of 
   some old walls of $F$.
   Thus one can find a family $\mathcal{G}$ of facets
   of $\DD$ such that $\sigma_p = \cap_{F' \in \mathcal{G}} (F' \cap F)$ and
   $F' \cap F \lessdot F$ for all $F' \in \mathcal{G}$.
   Since $|F' \cup F| = |F'| + 1 < n$, one knows that $F' \cap F$ is a facet of $\DD_F$
   of type (\ref{intersection}) for all $F' \in \mathcal{G}$.
   Let $\widetilde{F} = F_k \cap F$ be the last facet in the family 
   $\{ F' \cap F | F' \in \mathcal{G} \}$ (pick $k$ as small as possible).
   Since all facets of $\DD_F$ occurring earlier than $\widetilde{F}$ 
   have the form $F \cap F_i$ such that $F_i < F_k$ and $F_i \cap F \lessdot F$, 
   one can see
   $\widetilde{G} \subseteq \tilde{\sigma}_p \subseteq \widetilde{F}$.
   Thus $\widetilde{D}$
   is a shelling-trapped decomposition of $F \cup \{ v \}$ over $\DD_F$.
   Also one can define $\widetilde{w} \in \mathfrak{S}_{p-1}$ by
   $$
   \widetilde{w} (j) = \left\{
   \begin{array}{cl}
      w (j-1) & \text{if } 1< j \le p-1,\\
      w (p-1) & \text{if } j = 1.
   \end{array}
   \right.
   $$

   \vspace{2mm}
   \emph{Case} 2. $\sigma_p = F$.
   \vspace{2mm}
   
   In this case, we claim that
   $$
   \widetilde{D} = \{ ( \bar{\sigma}_1, \widetilde{F}_{i_1}), \dots,
   (\bar{\sigma}_k \cup \{ v \}, \widetilde{F}_{i_k}), \dots,
   ( \bar{\sigma}_{p-1}, \widetilde{F}_{i_{p-1}}) \}
   $$
   is a shelling-trapped decomposition of $F \cup \{ v \}$.

   Let $k = w (1)$. Define
   $$
   \tilde{\sigma}_j = \left\{ 
   \begin{array}{cl}
      \sigma_j \cap F & \text{ for } j = k, \\
      (\sigma_j \cap F) \cup \{ v \} & \text{ for } j = 1, \dots, \hat{k}, \dots, p-1 .
   \end{array}
   \right.
   $$
   Then
   $$
   \widetilde{G}_{i_k} 
   \subseteq \tilde{\sigma}_k = \sigma_k \cap F
   \subseteq \widetilde{F}_{i_k}
   $$
   and
   $$
   \widetilde{G}_{i_j} 
   \subseteq \tilde{\sigma}_j = (\sigma_j \cap F) \cup \{ v \} 
   \subseteq \widetilde{F}_{i_j}
   $$
   for $j = 1, \dots, \hat{k}, \dots, p-1$.
   Thus $\widetilde{D}$ is a shelling-trapped decomposition of $F \cup \{ v \}$.
   Also one can define $\widetilde{w} \in \mathfrak{S}_{p-2}$ as follows:
   $$
   \widetilde{w} (j) = \left\{
   \begin{array}{cl}
      w (j+1)       & \text{if } w (j+1) < k, \\
      w (j+1) -1   & \text{if } w (j+1) > k.
   \end{array}
   \right.
   $$

   Conversely, let 
   $$
   \widetilde{D} = \{ ( [F \cup \{ v \}] - \tilde{\sigma}_1, \widetilde{F}_{i_1}),
   \dots, ( [F \cup \{ v \}] - \tilde{\sigma}_p, \widetilde{F}_{i_p}) \}
   $$
   be a shelling-trapped decomposition of $F \cup \{ v \}$, where
   $\widetilde{F}_{i_1} < \dots < \widetilde{F}_{i_p}$ are facets of 
   $\DD_F$, and let $\widetilde{w}$ be a permutation in $\mathfrak{S}_{p-1}$.
   There is at most one facet of $\DD_F$ of type (\ref{intersection})
   because $\widetilde{F}_{i_j} \cup \widetilde{F}_{i_k} = F \cup \{ v \}$ for all $j \ne k$.
   Since $\widetilde{F}_{i_1} < \dots < \widetilde{F}_{i_p}$ and the facets of type (\ref{intersection})
   appear earlier than the ones of type (\ref{union}), there are two possible cases.
   
   \vspace{2mm}
   \emph{Case} 1. $v \notin \widetilde{F}_{i_1}$ and $v \in \widetilde{F}_{i_j}$
   for $j = 2, \dots, p$.
   \vspace{2mm}
   
   In this case, $v \notin \tilde{\sigma}_p$, i.e., $v \in [F \cup \{ v \}] - \tilde{\sigma}_p$.
   One can show that a family
   $$
   D = \{ ( [F \cup \{ v \}] - \tilde{\sigma}_2, F_{i_2}),
   \dots, ( [F \cup \{ v \}] - \tilde{\sigma}_p, F_{i_p}),
   ( [n] - \tilde{\sigma}_1, F) \},
   $$
   where $F_{i_j} = (\widetilde{F}_{i_j} \cup F) - \{ v \}$ for $j =2, \dots, p$,
   is a shelling-trapped decomposition of $[n]$ and $w \in \mathfrak{S}_{p-1}$
   is defined by
   $$
   w (j) = \left\{
   \begin{array}{cl}
      \widetilde{w} (j+1) & \text{if } 1 \le j < p-1 ,\\
      \widetilde{w} (1)    & \text{if } j = p-1.
   \end{array}
   \right.
   $$
   
   \vspace{2mm}
   \emph{Case} 2. $v \in \widetilde{F}_{i_j}$ for $j = 1, \dots, p$.
   \vspace{2mm}
   
   In this case, there is a $k$ such that $v \in [F \cup \{ v \}] - \tilde{\sigma}_k$.
   One can show that the family
   $$
   D = \{ ( F - \tilde{\sigma}_1, F_{i_1}),
   \dots, ( F - \tilde{\sigma}_p, F_{i_p}),
   ([n] - F, F) \},
   $$
   where $F_{i_j} = (\widetilde{F}_{i_j} \cup F) - \{ v \}$ for $j =2, \dots, p$,
   is a shelling-trapped decomposition of $[n]$ and $w \in \mathfrak{S}_p$
   can be defined by
   $$
   w (j) = \left\{
   \begin{array}{cl}
      \widetilde{w} (j-1)      & \text{if } 1 < j \text{ and } \widetilde{w} (j-1) < k, \\
      \widetilde{w} (j-1) +1  & \text{if } 1 < j \text{ and } \widetilde{w} (j-1) \ge k, \\
      k                & \text{if } j = 1.
   \end{array}
   \right.
   $$
   
   For the second assertion, one can show that each subspace in $\overline{C}_{D, w}$ is not changed under the map
   $\psi_1$ in the proof of Lemma~\ref{lemma-upper-interval} since it is the intersection of all subspaces
   that lie weakly below it and cover $U_{\overline{F}}$.
   We will show that $\psi(\overline{C}_{D, w} - U_{\overline{F}}) = \psi_2 (\overline{C}_{D, w} - U_{\overline{F}})$
   is a chain satisfying all conditions for $\overline{C}_{\widetilde{D}, \widetilde{w}}$.
   If $\bar{\tau}_j$ is a set satisfying $\overline{F}_{i_j} \subseteq \overline{\tau}_j \subseteq \bar{\sigma}_j$, then
   $$
   \psi (U_{\bar{\tau}_j}) = \left\{
   \begin{array}{cl}
      U_{\bar{\tau}_j}                       &\text{if } j = 1, 2, \dots, p-1,\\
      U_{[\bar{\tau}_p \cap F] \cup \{ v \}} &\text{if } j = p \text{ and } \bar{\tau}_p \ne \overline{F},\\
      \reals^{|F|+1}                         &\text{if } j = p \text{ and } \bar{\tau}_p = \overline{F},
   \end{array}
   \right.
   $$
   where the coordinates of $\reals^{|F|+1}$ are indexed by $F \cup \{ v \}$.

   There are two cases:

   \vspace{2mm}
   \emph{Case} 1. $\sigma_p \ne F$.
   \vspace{2mm}
   
   In this case, one can see that $F - \widetilde{F}$ has only one element, say $x$, and $\bar{\sigma}_p$ contains $x$.
   Thus one can choose $U_{\overline{F}_{i_1} / \cdots / \overline{F}_{i_{p-1}} / \overline{F} \cup \{ x \}}$ 
   as the second subspace in the chain $\widehat{C}_w$.
   Then $\psi (U_{\overline{F}_{i_1} / \cdots / \overline{F}_{i_{p-1}} / \overline{F} \cup \{ x \}})
   = U_{[F \cup \{ v \}] - \widetilde{F} / [F \cup \{ v \}] - \widetilde{F}_{i_1} / \cdots / [F \cup \{ v \}] - \widetilde{F}_{i_{p-1}}}$
   since $\overline{F}_{i_j} = [F \cup \{ v \}] - \widetilde{F}_{i_j}$ for $j = 1, 2, \dots, p-1$ and 
   $[F \cup \{ v \}] - \widetilde{F} = \{ v, x \}$.
   Moreover, the image of $\widehat{C}_w - U_{\overline{F}_{i_1}/\cdots /\overline{F}_{i_p}}$ under $\psi$ is a saturated chain in 
   $[U_{[F \cup \{ v \}] - \widetilde{F} / [F \cup \{ v \}] - \widetilde{F}_{i_1} / \cdots / [F \cup \{ v \}] - \widetilde{F}_{i_{p-1}}},
   U_{F \cup \{ v \}}]$ whose covering relations of type (a) (in Step (1) of Section~\ref{sec-chain}) 
   correspond to the covering relations in the chain $C_{\widetilde{w}}$ in $\Pi_p$.
   Therefore, the image of $\overline{C}_{D, w}$ under $\psi$ can be chosen as $\overline{C}_{\widetilde{D}, \widetilde{w}}$.

   \vspace{2mm}
   \emph{Case} 2. $\sigma_p = F$.
   \vspace{2mm}
   
   In this case, it is not hard to see that the image of $\widehat{C}_w$ under $\psi$ is the saturated chain in 
   $[U_{[F \cup \{ v \}] - \widetilde{F}_{i_1} / \cdots / [F \cup \{ v \}] - \widetilde{F}_{i_{p-1}}},
   U_{F \cup \{ v \}}]$
   whose covering relations of type (a) (in Step (1) of Section~\ref{sec-chain})
   correspond to the covering relations in the chain $C_{\widetilde{w}}$ in $\Pi_{p-1}$.
   Thus the image of $\overline{C}_{D, w}$ under the map $\psi$ can be chosen as $\overline{C}_{\widetilde{D}, \widetilde{w}}$.
\end{proof}

\begin{example}
   Let $\DD$ be the shellable simplicial complex from Example~\ref{ex-upper-iso-last}.
   In Example~\ref{ex-upper-iso-chain}, $C_{D, w}$ is the chain 
   $$
   \hat{0} \lessdot U_{67} \lessdot U_{23 / 67} \lessdot
   U_{45 / 23 / 67}
   \lessdot U_{45 / 123 / 67} \lessdot
   U_{45 / 12367} \lessdot U_{1234567}
   $$
   for the shelling-trapped decomposition
   $$
   D = \{ (45, F_1 = 12367), (123, F_6 = 14567), (67, F_7 = 12345) \}
   $$
   of $[7]$ and the permutation $w$ in $\mathfrak{S}_2$ 
   with $w (1) = 2$ and $w (2) = 1$.
   
   Since $67 = \overline{F}_7$, the corresponding shelling-trapped decomposition 
   $\widetilde{D}$ of the set $\{ 1, 2, 3, 4, 5, v \}$ is
   $$
   \widetilde{D} = \{ (45, \widetilde{F}_1 = 123v), 
   (123v, \widetilde{F}_6 = 145v) \}
   $$
   and the corresponding permutation $\widetilde{w} \in \mathfrak{S}_1$ is the identity.
   
   The map $\psi$ from the proof of Lemma~\ref{lemma-upper-interval} sends the chain 
   $$
   U_{23 / 67} \lessdot
   U_{45 / 23 / 67}
   \lessdot U_{45 / 123 / 67} \lessdot
   U_{45 / 12367}
   $$
   to the chain
   $$
   U_{23} \lessdot U_{45 / 23}
   \lessdot U_{45 / 123} \lessdot U_{45 / 123v}.
   $$
   and this chain satisfies the conditions for $\overline{C}_{\widetilde{D}, \widetilde{w}}$.

   The intersection lattice for $\DD_F$ is shown in Figure~\ref{pic-iso-upper-interval}
   and the chain $\overline{C}_{\widetilde{D}, \widetilde{w}}$ is represented by thick
   lines.
   
   \begin{figure}
   \begin{center}
   \includegraphics[width=0.9\textwidth]{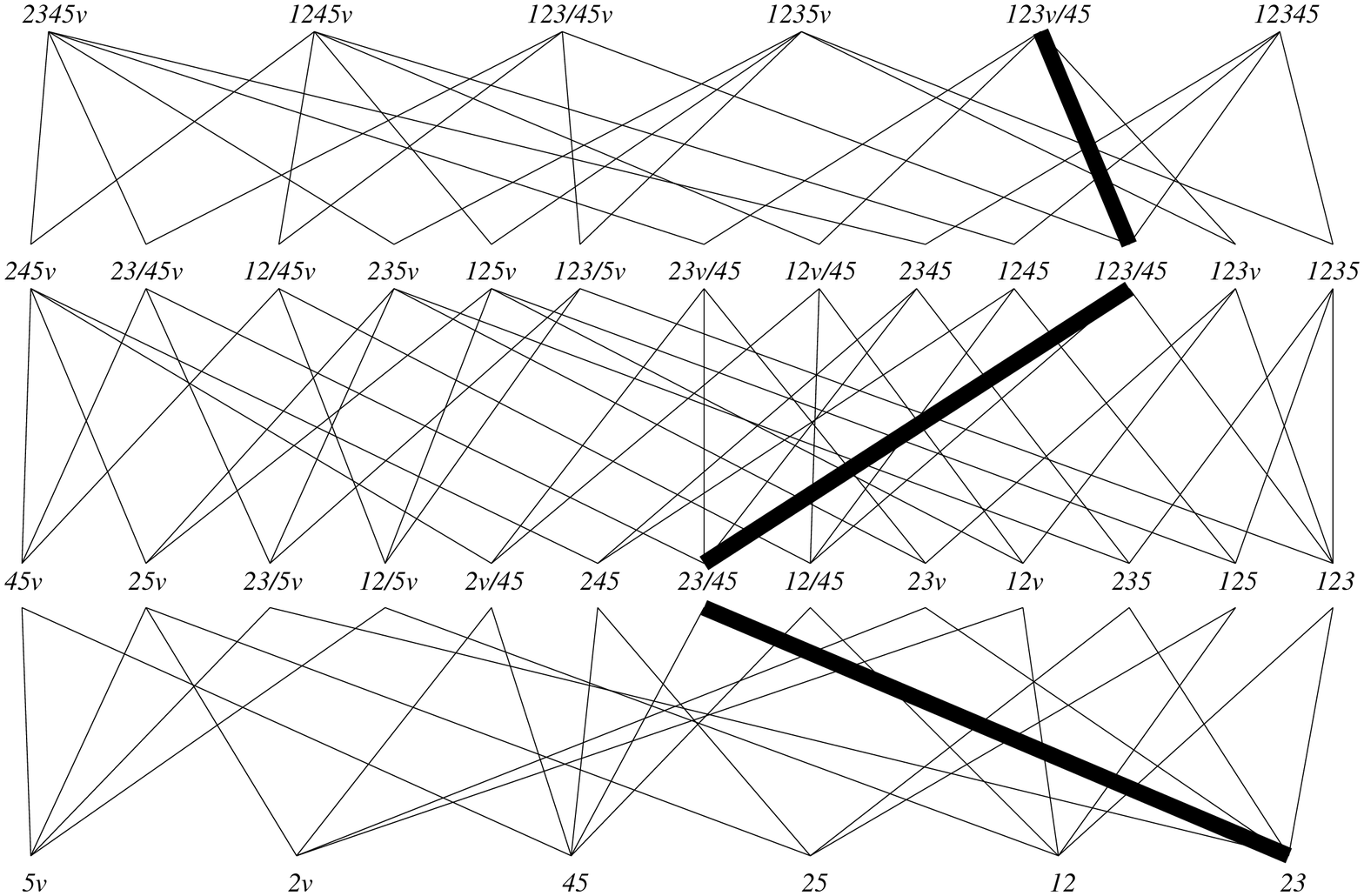}
   \begin{caption}
      {The interval $(\hat{0}, \hat{1})$ in $L_{\DD_F}$}
      \label{pic-iso-upper-interval}
   \end{caption}
   \end{center}
   \end{figure}
\end{example}

\emph{Proof of Theorem \ref{one-set-nonpure}. }
   By Lemma~\ref{iso-to-links}, it is enough to show the assertion for the case when
   $\bar{\sigma} = \cup_{i=1}^q \overline{F}_i = [n]$.
   Since every chain $\overline{C}_{D, w}$ is saturated, it is enough to show that
   $\WHD (\overline{L}_\DD)$, the simplicial complex obtained after removing
   the corresponding simplices for all
   pairs $(D, w)$, is contractible.
   We use induction on the number $q$ of facets of $\DD$.

   \emph{Base case}: $q=2$.
   If $\DD$ has only two facets $F_1$ and $F_2$ and
   $\overline{F}_1 \cup \overline{F}_2 = [n]$, 
   then $F_2$ has only one element and $G_2 = \emptyset$.
   It is easy to see that the order complex $\DD(\overline{L}_\DD)$ 
   is homotopy equivalent to $\mathbb{S}^0$ and
   $D = \{ ( [n], F_2) \}$ is the only 
   shelling-trapped decomposition of $[n]$,
   while $\overline{C}_{D, \emptyset} = (U_{\overline{F}_2})$ 
   is the corresponding saturated chain.
   Therefore, $\WHD (\overline{L}_\DD)$ is contractible when $q=2$.
   
   \emph{Inductive step}.
   Now, assume that $\WHD (\overline{L}_\DD)$ is contractible 
   for all shellable simplicial complexes $\DD$ with less than $q$ facets.
   For simplicity, denote $L = L_{\DD}$.
   Let $F = F_q$ be the last facet in the shelling order of $\DD$
   and $H = U_{\overline{F}}$.
   Let $L'$ be the intersection lattice for $\DD '$,
   where $\DD '$ is the simplicial complex generated by the facets
   $F_1, \dots, F_{q-1}$.
   
   Let $\overline{L}_{\ge H}$ denote the subposet of elements in $\overline{L}$
   which lie weakly above $H$.
   Consider the decomposition of $\WHD (\overline{L}) = X \cup Y$,
   where 
   $X$
   is the simplicial complex 
   obtained by removing all simplices corresponding to chains  
   $\overline{C}_{D, w}$ and $\overline{C}_{D, w} - H$ 
   from $\DD(\overline{L}_{\ge H})$
   for all $\overline{C}_{D, w}$ containing $H$, and
   $Y$
   is the simplicial complex
   obtained by removing all simplices corresponding to chains
   $\overline{C}_{D, w}$ not containing $H$ from $\DD(\overline{L} - \{ H \})$.
   Our goal will be to show that $X, Y$ and $X \cap Y$ are all contractible,
   and hence so is $X \cup Y (= \WHD (L))$.

   \vspace{2mm}
   \emph{Step 1. Contractibility of $X$}
   \vspace{2mm}
   
   Since $X$ has a cone point $H$, it is contractible.
   
   \vspace{2mm}
   \emph{ Step 2. Contractibility of $Y$}
   \vspace{2mm}
   
   Define the closure relation $\pi$ on $L$ which sends a subspace
   to the join of the elements covering $\hat{0}$ which lie below it except $H$.
   Then the closed sets form a sublattice of $L$, which is the intersection lattice
   $L'$ for the diagonal arrangement corresponding to $\DD '$. 
   It is well-known that the inclusion of closed sets
   $\overline{L} \cap L' \hookrightarrow \overline{L} - \{ H \}$
   is a homotopy equivalence (see~\cite[Lemma 7.6]{BjornerWachs1}).
   We have to consider the following two cases:
   
   \vspace{2mm}
   \emph{Case} 1.    
   $S := \cup_{i=1}^{q-1} \overline{F}_i \ne [n]$.
   \vspace{2mm}
   
   Then $\overline{L} \cap L' = L' - \{ \hat{0} \}$ since $U_S \in \overline{L}$.
   Since $\overline{L} \cap L'$ has a cone point $U_S$, it is contractible.
   Since $S \ne [n]$, there is no shelling-trapped decomposition of $[n]$ for $\DD '$.
   Thus $Y$ is contractible.
      
   \vspace{2mm}
   \emph{Case} 2.
   $S := \cup_{i=1}^{q-1} \overline{F}_i  = [n]$
   \vspace{2mm}
   
   In this case, $\overline{L} \cap L' = \overline{L'}$.
   Moreover, $\WHD (\overline{L'})$ is homotopy equivalent to
   $Y$ since every element in a chain 
   $\overline{C}_{D, w}$ in $\WHD (\overline{L} - \{ H \})$ is fixed under $\pi$.
   Since $\DD'$ has $q-1$ facets, the induction hypothesis implies that
   $\WHD (\overline{L'})$ is contractible and hence so is $Y$.
   
   \vspace{2mm}
   \emph{Step 3. Contractibility of $X \cap Y$}
   \vspace{2mm}
   
   Note that $X \cap Y$ is obtained by removing simplices corresponding to
   $\overline{C}_{D, w} - H$ for all $\overline{C}_{D, w}$ containing $H$
   from $\overline{L}_{>H}$.
   By Lemma \ref{upper-last-facet}, $(\overline{L})_{> H}$ is 
   isomorphic to the proper part of the intersection lattice $L_F$ for 
   the diagonal arrangement corresponding to $\DD_F$ on $F \cup \{ v \}$.
   Also, Lemma~\ref{STD-correspondence} implies that 
   $X \cap Y$ is isomorphic to 
   $\WHD (\overline{L}_F)$, where $\WHD (\overline{L}_F)$ is obtained by 
   removing simplices corresponding to 
   $\overline{C}_{\widetilde{D}, \widetilde{w}}$ for all 
   shelling-trapped decomposition $\widetilde{D}$ of $F \cup \{ v \}$ 
   and $\widetilde{w} \in \mathfrak{S}_{|\widetilde{D}|-1}$ from $\overline{L}_F$.
   Since $\DD_F$ has fewer facets than $\DD$, the induction hypothesis 
   implies $\WHD (\overline{L}_F)$ is contractible and hence 
   $X \cap Y$ is also contractible.
\qed

\begin{example}
\label{ex-shelling-trapped-decomposition}
   Let $\DD$ be a simplicial complex from Example~\ref{ex-upper-iso}.
   Then one can see that
   $$
   \begin{array}{ccccc}
      F_1=123, &F_2=234, &F_3=35, &F_4=45
   \end{array}
   $$
   is a shelling and 
   $$
   \begin{array}{ccccc}
      G_1 = 123, &G_2 = 23, &G_3 = 3, &G_4 = \emptyset.
   \end{array}
   $$
   Let $\bar{\sigma} = 12345$.
   Then there are two possible (unordered) shelling-trapped decompositions 
   of $\bar{\sigma}$ 
   (see Table \ref{delta-table}).

   Thus, Theorem~\ref{one-set-nonpure} implies 
   $\DD(\hat{0}, U_{12345})$ is homotopy equivalent to a wedge of 
   two circles.
   The intersection lattice $L_\DD$ and the order complex for its proper part are shown in 
   Figure~\ref{intersection-lattice}.
   Note that the chains $\overline{C}_{D, w}$ and the simplices corresponding to each 
   shelling-trapped decomposition are represented by thick lines.

   \begin{figure}[!t]
   \begin{center}
      \begin{tabular}{ccc}
         \includegraphics[width=0.5\textwidth]{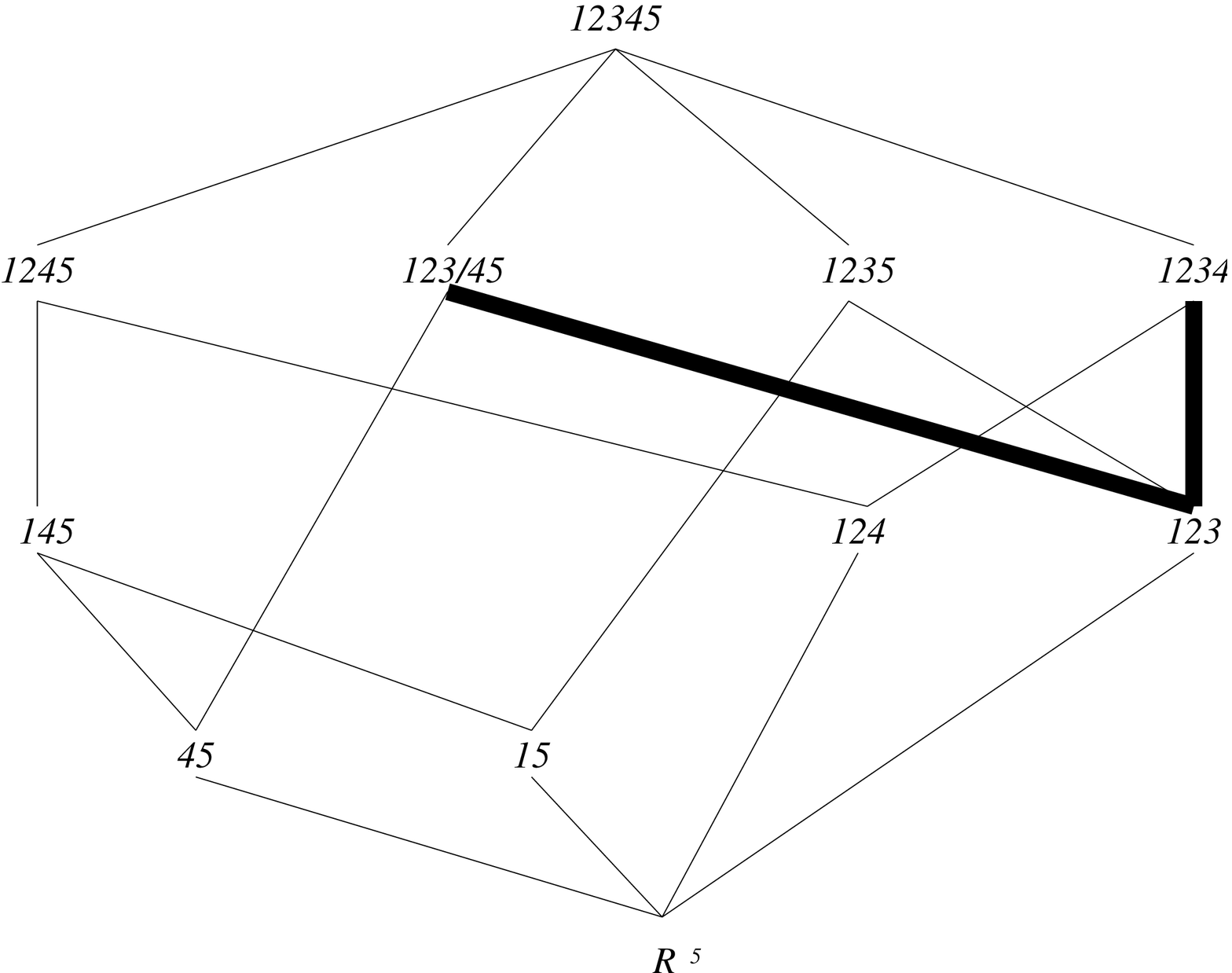}& &
         \includegraphics[width=0.35\textwidth]{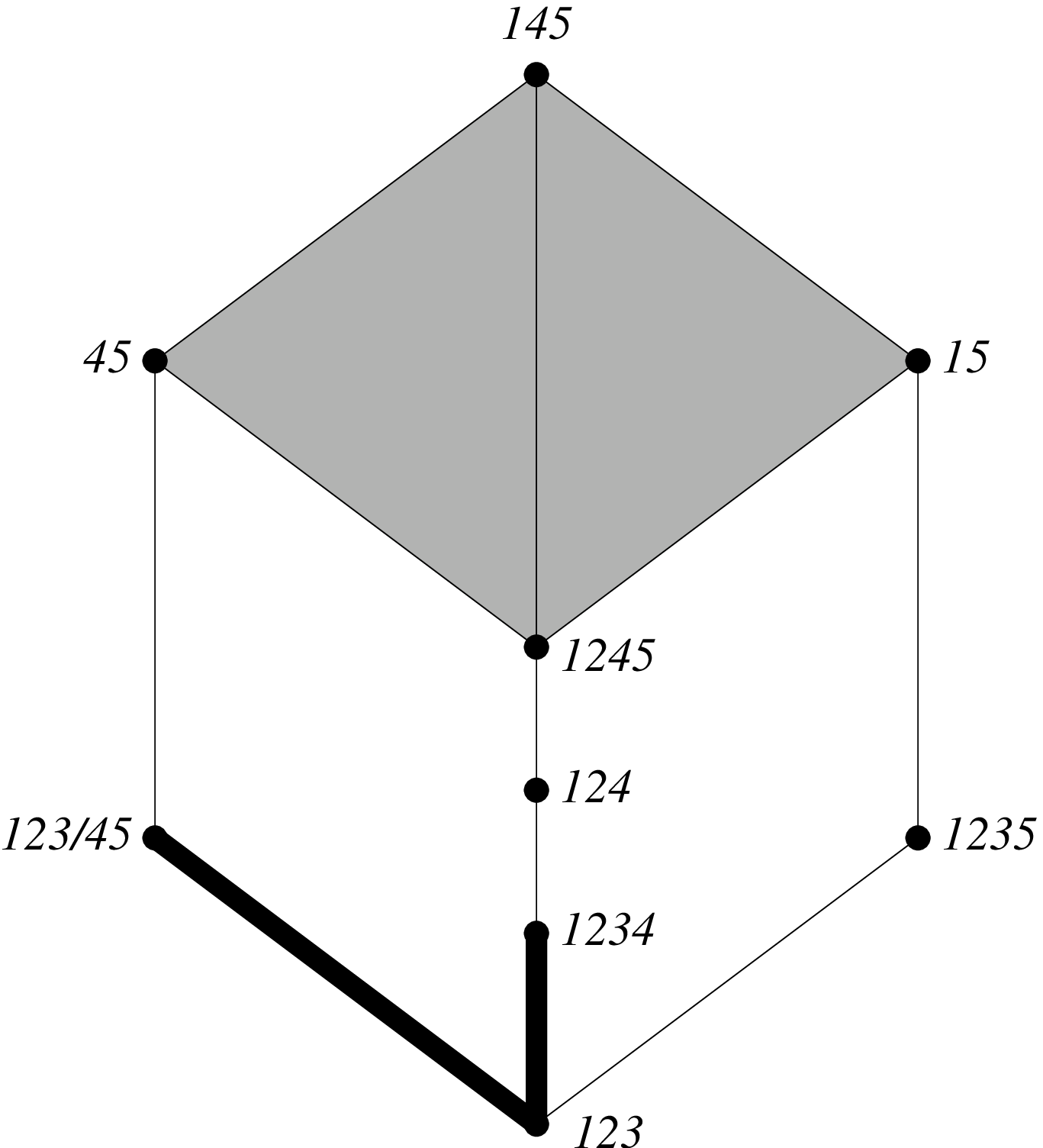}\\
         {\small (a) $L_\DD$} & \hspace{0.5cm} & {\small (b) The order complex of $\overline{L}_\DD$} \\
      \end{tabular}

      \begin{caption}
         {The intersection lattice $L_\DD$ of $\mathcal{A}_\DD$ and the order complex of its proper part}
         \label{intersection-lattice}
      \end{caption}
   \end{center}
   \end{figure}
\end{example}

\section{The homotopy type of the singularity link of $\mathcal{A}_\DD$}
\label{set-homotopy}

In this section, we give the corollary about the homotopy type of the singularity link of 
$\mathcal{A}_\DD$ when $\DD$ is shellable. 
We also give the homology version of the corollary.

Ziegler and \v{Z}ivaljevi\'{c}~\cite {ZieglerZivaljevic} showed the following theorem
about the homotopy type of $\mathcal{V}_\mathcal{A}^\circ$.

\begin{theorem}
\label{Ziegler-Zivaljevic-theorem}
   For every subspace arrangement $\mathcal{A}$ in $\reals^n$,
   $$
   \mathcal{V}_\mathcal{A}^\circ \simeq
   \bigvee_{x \in L_{\mathcal{A}} - \{ \hat{0} \}}
   ( \DD ( \hat{0}, x ) * \mathbb{S}^{\dim (x) - 1} ).
   $$
\end{theorem}

From this and our results in Section~\ref{sec-homotopy}, one can deduce
the following.

\begin{corollary}
\label{cor-singularity-link}
   Let $\DD$ be a shellable simplicial complex on $[n]$.
   The singularity link of $\mathcal{A}_{\DD}$ has the homotopy type of a wedge of
   spheres, consisting of $p!$ spheres of dimension 
   $$
   n+p(2-n)+\sum_{j=1}^p |F_{i_j}|-2
   $$
   for each shelling-trapped decomposition
   $\{ (\bar{\sigma}_1, F_{i_1}), \ldots, (\bar{\sigma}_p, F_{i_p}) \}$
   of some subset of $[n]$.
\end{corollary}

\noindent
\emph{Proof sketch.}
This is a straightforward, but tedious, calculation.
By Theorem~\ref{Ziegler-Zivaljevic-theorem}, one needs to understand homotopy types
of $\DD(\hat{0}, H)$ for $H \in L_\DD$.
Lemmas~\ref{element-of-lattice} and \ref{more-than-one} reduce this to the case of
$\DD(\hat{0}, U_{\bar{\sigma}})$, which is described fully by
Theorem~\ref{one-set-nonpure}.
The rest is some bookkeeping about shelling-trapped decompositions.
\qed

\begin{table}[!t]
   \begin{center}
   \begin{tabular}{|c|c||c|c|}
      \hline
      shelling-trapped decomp. & $\dim$ & shelling-trapped decomp.   & $\dim$ \\ \hline \hline
      $\{(45, F_1)\}$          & $3$    & $\{(15, F_2)\}$            & $3$    \\ \hline
      $\{(145, F_2)\}$         & $3$    & $\{(124, F_3)\}$           & $2$    \\ \hline
      $\{(1245, F_3)\}$        & $2$    & $\{(123, F_4)\}$           & $2$    \\ \hline
      $\{(1234, F_4)\}$        & $2$    & $\{(1235, F_4)\}$          & $2$    \\ \hline
      $\{(12345, F_4)\}$       & $2$    & $\{(45,F_1), (123, F_4)\}$ & $2$    \\ \hline
   \end{tabular}
   \end{center}
      \begin{caption}
      {Shelling-trapped decompositions for $\DD$}
      \label{delta-table}
   \end{caption}
\end{table}

\begin{example}
\label{ex-homology-of-link}
   Let $\DD$ be a simplicial complex from Example~\ref{ex-upper-iso}.
   In Example~\ref{ex-shelling-trapped-decomposition}, we show
   $$
   F_1=123, F_2=234, F_3=35, F_4=45
   $$ 
   is a shelling of $\DD$ and 
   $$
   G_1 = 123, G_2 = 23, G_3 = 3, G_4 = \emptyset.
   $$
   Table~\ref {delta-table} shows shelling-trapped decompositions
   $\{(\bar{\sigma}_1, F_{i_1}), \ldots , (\bar{\sigma}_p, F_{i_p})\}$
   of subsets of $[5]$ with corresponding dimensions
   $$
   n + p(2-n) + \sum_{j=1}^p |F_{i_j}|-2.
   $$
   Therefore Corollary~\ref{cor-singularity-link} shows that
   the singularity link of $\mathcal{A}_\DD$ is homotopy equivalent to
   a wedge of three $3$-dimensional spheres and eight $2$-dimensional spheres.
\end{example}

The following theorem is a homology version of Corollary~\ref{cor-singularity-link}.

\begin{theorem}
\label{nonpure-by-walls}
   Let $\DD$ be a shellable simplicial complex and
   $F_1, \ldots, F_q$ be a shelling of $\DD$.
   Then
   $\dim_\FF H_i (\mathcal{V}^\circ_{\mathcal{A}_{\DD}} ; \FF)$
   is the number of ordered shelling-trapped decompositions
   $((\bar{\sigma}_1, F_{i_1}), \ldots , (\bar{\sigma}_p, F_{i_p}))$
   with
   $$
   i = n + p(2-n) + \sum_{j=1}^p |F_{i_j}|-2.
   $$
\end{theorem}

This can be proven without Theorem~\ref{one-set-nonpure} by combining a result of
Peeva, Reiner and Welker~\cite[Theorem 1.3] {PeevaReinerWelker} with results of 
Herzog, Reiner and Welker~\cite[Theorem 4, Theorem 9] {HerzogReinerWelker99} 
along with the theory of \emph{Golod rings}.
It is what motivated us to prove the stronger Corollary~\ref{cor-singularity-link} 
and eventually Theorem~\ref{main-theorem}.

\section{$K(\pi, 1)$ examples from matroids}
\label{sec-examples}

In this section, we give some examples of diagonal arrangements $\mathcal{A}$
where the complement $\mathcal{M}_\mathcal{A}$ is $K(\pi,1)$, 
coming from rank $3$ matroids.

One should note that an arrangement having any subspace of real codimension $1$
(hyperplane) will have $\mathcal{M}_\mathcal{A}$ disconnected.
So one may assume without loss of generality that all subspaces have
real codimension at least $2$.
Furthermore, if any maximal subspace $U$ in $\mathcal{A}$ has codimension
at least $3$, then it is not hard to see that $\mathcal{M}_\mathcal{A}$
is not $K(\pi,1)$.
Hence we may assume without loss of generality that all maximal subspaces
have real codimension $2$.

A hyperplane arrangement $\mathcal{H}$ in $\reals^n$ is \emph{simplicial} if
every chamber in $M_\mathcal{H}$ is a simplicial cone.
Davis, Januszkiewicz and Scott \cite{DavisJanuszkiewiczScott} showed the following theorem.

\begin{theorem}
\label{Davis-Januszkiewicz-Scott-Thm}
   Let $\mathcal{H}$ be a simplicial real hyperplane arrangement in $\reals^n$.
   Let $\mathcal{A}$ be any arrangement of codimension-$2$ subspaces in $\mathcal{H}$ 
   which intersects every chamber in a codimension-$2$ subcomplex.
   Then $\mathcal{M}_{\mathcal{A}}$ is $K(\pi, 1)$.
\end{theorem}

\begin{remark}
   In order to apply this to diagonal arrangements, 
   we need to consider hyperplane arrangements $\mathcal{H}$ which are
   subarrangements of the real braid arrangement $\mathcal{B}_n$ in $\reals^n$ and also simplicial.
   It turns out (and we omit the straightforward proof) that all such
   arrangements $\mathcal{H}$ are direct sums of smaller braid arrangements.
   So we only consider $\mathcal{H} = \mathcal{B}_n$ itself here.
\end{remark}

\begin{corollary}
\label{Davis-Januszkiewicz-Scott-Cor}
   Let $\mathcal{A}$ be a subarrangement of the $3$-equal arrangement of $\reals^n$
   so that
   $$
   \mathcal{A} = \{ U_{i j k} \mid \{ i, j, k \} \in T_{\mathcal{A}} \},
   $$
   for some collection $T_{\mathcal{A}}$ of $3$-element subsets of $[n]$.
   Then $\mathcal{A}$ satisfies the hypothesis of 
   Theorem~\ref{Davis-Januszkiewicz-Scott-Thm} 
   (and hence $\mathcal{M}_{\mathcal{A}}$ is $K(\pi, 1)$) if and only if
   every permutation $w = w_1 w_2 \dots w_n$ in $\mathfrak{S}_n$
   has at least one triple in $T_{\mathcal{A}}$ consecutive, i.e.,
   there exists $j$ such that $\{ w_j, w_{j+1}, w_{j+2} \} \in T_{\mathcal{A}}$.
\end{corollary}

\begin{proof}
   It is easy to see that there is a bijection between chambers of 
   the real Braid arrangement $\mathcal{B}_n$ in $\reals^n$ and
   permutations $w = w_1 \cdots w_n$ in $\mathfrak{S}_n$.
   Moreover, each chamber has the form $x_{w_1} > \cdots > x_{w_n}$
   with bounding hyperplanes 
   $$
   x_{w_1} = x_{w_2}, x_{w_2} = x_{w_3},\dots, x_{w_{n-1}} = x_{w_n}
   $$
   and intersects the $3$-equal subspaces of the form
   $x_{w_i} = x_{w_{i+1}} = x_{w_{i+2}}$ for $i = 1, 2, \dots, n-2$.
\end{proof}

We seek shellable simplicial complexes $\DD$ for which 
$\mathcal{A}_\DD$ satisfies this condition.

If $\DD$ is the independent set complex $\mathcal{I}(M)$ for some matroid $M$
(see \cite{Oxley} for the definition of independent sets and further background on matroids), 
then facets of $\DD$ are bases of $M$.
Simplicial complexes of this kind are called \emph{matroid complexes},
and they are known to be shellable~\cite{Bjorner92}.
For a rank $3$ matroid $M$ on $[n]$, consider 
$$
\begin{aligned}
\mathcal{A}_{\mathcal{I}(M^\perp)} 
&= \{ U_{ijk} \mid \{ i, j, k \} = [n] - B \text{ for some } B \in \mathcal{B}(M^\perp) \}\\
&= \{  U_{ijk} \mid \{ i, j, k \} \in \mathcal{B}(M) \},
\end{aligned}
$$
where $M^\perp$ is the dual matroid of $M$.
Note that adding a loop to $M$ does not change the structure of
the intersection lattice for $\mathcal{A}_{\mathcal{I}(M^\perp)}$.
Thus, if $\mathcal{B}(M)$ on the set of all non-loop elements satisfies
the condition of Corollary~\ref{Davis-Januszkiewicz-Scott-Cor}, then
the diagonal arrangement corresponding to the matroid 
in which all loops have been deleted has $K(\pi,1)$ complement,
and hence $\mathcal{A}_{\mathcal{I}(M^\perp)}$ has $K(\pi,1)$ complement.

\begin{definition}
   Let $M$ be a rank $3$ matroid on $[n]$.
   Say $M$ is \emph{DJS} if $\mathcal{B}(M)$ on the set of all non-loop elements 
   satisfies the condition of Corollary~\ref{Davis-Januszkiewicz-Scott-Cor}.
\end{definition}

The following example shows that matroid complexes are not DJS in general.
Thus we look for some subclasses of matroid complexes which are DJS, and hence
whose corresponding diagonal arrangements have $K(\pi, 1)$ complements;
for these, Theorem~\ref{nonpure-by-walls} gives us the group cohomology 
$H^\bullet (\pi, \mathbb{Z})$.

\begin{example}
   Let $M$ be a matroid on $[6]$
   which has three distinct parallel classes $\{ 1, 6 \},$ $\{ 2, 4 \}$ and $\{ 3, 5 \}$.
   Then $M$ is self-dual and $\mathcal{I}(M^\perp)$ is a simplicial complex on $[6]$ 
   whose facets are $123, 134, 145,$ $125, 236, 256, 346$ and $456$.
   But $w = 124356$ is a permutation that does not satisfy the condition of
   Corollary~\ref{Davis-Januszkiewicz-Scott-Cor}.
\end{example}

Recall that a matroid is \emph{simple} if it has no loops nor parallel elements.
The following proposition shows that rank $3$ simple matroids are DJS.

\begin{proposition}
   Let $M$ be a matroid of rank $3$ on $[n]$. 
   If $M$ does not have parallel elements, then $M$ is DJS. 
   In particular, rank $3$ simple matroids are DJS.
\end{proposition}

\begin{proof}
   Without loss of generality, we may assume that $M$ is simple.
   $M$ is not DJS if and only if there is a permutation 
   $w \in \mathfrak{S}_n$ such that every consecutive triple is not in $\mathcal{B}(M)$.
   Since $M$ is simple,
   the latter statement is true if and only if each consecutive triple in $w$ 
   forms a circuit, i.e., all elements lie on a rank $2$ flat.
   But this is impossible since $M$ has rank $3$.
\end{proof}

The following two propositions give some subclasses of matroids with parallel elements 
which are DJS.

\begin{proposition}
   Let $M$ be a rank $3$ matroid on $[n]$ with no circuits of size $3$.
   Let $P_1, \dots, P_k$ be the distinct parallel classes which have more than 
   one element, and let $N$ be the set of all non-loop elements which are not parallel 
   with anything else.
   Then, $M$ is DJS if and only if 
   $$
   \lfloor \frac{|P_1|}{2} \rfloor + \cdots + \lfloor \frac{|P_k|}{2} \rfloor - k < |N|-2.
   $$
\end{proposition}

\begin{proof}
   We may assume that $M$ does not have loops.
   Since $M$ does not have loops nor circuits of size $3$,
   $M$ is not DJS if and only if one can construct a permutation $w \in \mathfrak{S}_n$
   such that for each consecutive triple in $w$ 
   there are at least two elements which are parallel.
   This means if $w_i \in N$, then $w_{i-2}, w_{i-1}, w_{i+1}$ and $w_{i+2}$ 
   (if they exist) must be in the same parallel class. 
   Such a $w$ can be constructed if and only if 
   $\lfloor \frac{|P_1|}{2} \rfloor + \cdots + \lfloor \frac{|P_k|}{2} \rfloor - k \ge |N|-2$.
\end{proof}

A simplicial complex $\DD$ on $[n]$ is \emph{shifted} if, for any face $\sigma$ of $\DD$, 
replacing any vertex $i \in \sigma$ by a vertex $j < i$ with $j \notin \sigma$
gives another face in $\DD$.
A matroid $M$ is \emph{shifted} if its independent set complex is shifted.
Klivans~\cite{Klivans} showed that a rank $3$ shifted matroid
on the ground set $[n]$ is indexed by some set $\{ a , b , c \}$ 
with $1 \le a < b < c \le n$ as follows:
$$
\mathcal{B}(M) = \{ (a', b', c') : 1 \le a' < b' < c' \le n, a' \le a, b' \le b, c' \le c \}.
$$

It is not hard to check the following.

\begin{proposition}
   Let $M$ be the shifted rank $3$ matroid on the ground set $[n]$ indexed by $\{ a, b, c \}$. 
   Then, $M$ is DJS if and only if $\lfloor \frac{c - b}{2} \rfloor < a$.
\end{proposition}

We have not yet been able to characterize all rank $3$ matroids which are DJS.

\section*{Acknowledgments}

The author would like to thank his advisor, Vic Reiner, for introducing the topic of subspace arrangements
and encouraging him to work on this problem.
The author also thanks Ezra Miller, Dennis Stanton, Michelle Wachs and Volkmar Welker for 
valuable suggestions.
Finally, the author thanks the referee for a careful reading of the manuscript, and valuable
suggestions for improvement.

\end{document}